\documentclass[a4paper,12pt]{article}

\pdfoutput=1

\usepackage[paper=letterpaper,margin=1in]{geometry}
\usepackage{graphics,epsfig, subfigure}
\usepackage{url,fullpage}
\usepackage{amssymb,amsmath,latexsym,amsfonts,amscd,amsthm}
\usepackage{latexsym,amsmath,amsfonts,amscd,amsthm,upref,graphics,graphicx,amssymb,hyperref}
\usepackage{fancyhdr}

\usepackage[numbers]{natbib}

\renewcommand{\bibpreamble}{NB. We place a superscript $^{>}$ for all references which have appeared after the date of the delivery of the lectures, as they are added in the editorial process.}

\newcommand{\nc}{\newcommand}
\newcommand{\nn}{\nonumber}

\newcommand{\comment}[1]{}

\nc{\G}{{\Gamma}}
\nc{\fH}{{\mathfrak{H}}}

\nc{\BC}{{\mathbb C}} \nc{\BQ}{{\mathbb Q}} \nc{\BR}{{\mathbb R}}
\nc{\BZ}{{\mathbb Z}} \nc{\BP}{{\mathbb P}} \nc{\BN}{{\mathbb N}}

\nc{\lcm}{\mathop{{\rm lcm}}}

\nc{\PS}{{\mbox{PSL}_2(\BZ)}} \nc{\SL}{{\mbox{SL}_2(\BZ)}}
\nc{\SR}{{\mbox{SL}_2(\BR)}} \nc{\PR}{{\mbox{PSL}_2(\BR)}}
\nc{\GL}{{\mbox{GL}_2^+(\BQ)}} \nc{\PQ}{{\mbox{PGL}_2^+(\BQ)}}
\nc{\GR}{{\mbox{GL}_2^+(\BR)}} \nc{\PG}{{\mbox{PGL}_2^+(\BR)}}
\nc{\GC}{{\mbox{GL}_2(\BC)}}

\newcommand{\tmat}[1]{{\tiny \left(\begin{matrix} #1 \end{matrix}\right)}}
\newcommand{\mat}[1]{\left(\begin{matrix} #1 \end{matrix}\right)}
\newcommand{\gen}[1]{\langle #1 \rangle}

\newcommand{\cO}{{\cal O}}
\newcommand{\IM}{\mathbb{M}}

\newcommand{\IQ}{\mathbb{Q}}

\newcommand{\IC}{\mathbb{C}}

\newcommand{\II}{\mathbb{I}}
\newcommand{\IZ}{\mathbb{Z}}

\newcommand{\im}{{\rm Im}}

\newtheorem{theorem}{\bf THEOREM}
\newtheorem{conjecture}{\bf CONJECTURE}

\renewenvironment{thebibliography}[1]{%
\begin{oldthebibliography}{#1}%
\setlength{\parskip}{0ex}%
\setlength{\itemsep}{0ex}%
}%
{%
\end{oldthebibliography}%
}

\begin{document}
~\\

\vskip 1cm

\begin{center}
{\large \bf Kashiwa Lectures on }\\~\\
{\Large \bf New Approaches to the Monster}
\end{center}
\medskip

\centerline{{\large John McKay}$^1$}
~\\
\centerline{edited and annotated by Yang-Hui He$^{2,3,4,5}$}

{\it {\small
\begin{tabular}{cl}
${}^{1}$ 
&CICMA \& Department of Mathematics and Statistics,
Concordia University, \\
&1455 de Maisonneuve Blvd.~West,
Montreal, Quebec, H3G 1M8, Canada\\
&\\
${}^{2}$ 
& London Institute for Mathematical Sciences, Royal Institution of Great Britain, \\
& 21 Albemarle Street, Mayfair, London W1S 4BS, UK; 
\\
${}^{3}$ & Merton College, University of Oxford, OX14JD, UK;
\\
${}^{4}$ & Department of Mathematics, City, University of London, EC1V 0HB, UK;
\\
${}^{5}$ & School of Physics, NanKai University, Tianjin, 300071, P.R.~China\\
\end{tabular}
}}

\begin{center}
mckay@encs.concordia.ca \quad hey@maths.ox.ac.uk
\end{center}

\vspace*{1.0ex}
\centerline{\textbf{Abstract}} \bigskip
These notes stem from lectures given by the first author (JM) at the 2008 
``Moonshine Conference in Kashiwa'' 
\footnote{Organized by the Institute for the Physics and Mathematics of the Universe (IPMU) under the support of the Graduate School of Mathematical Sciences, the University of Tokyo.}
and contain a number of new perspectives and observations on Monstrous Moonshine.
Because many new points have not appeared anywhere in print, it is thought expedient to update, annotate and clarify them (as footnotes), an editorial task which the second author (YHH) is more than delighted to undertake.
We hope the various puzzles and correspondences, delivered in a personal and casual manner,  will serve as
diversions intriguing to the community.

\newpage
\tableofcontents
\newpage


\section{Introduction}
I am very honoured to be able to attend and participate in this conference which I would very much appreciate to be a relevantly informal business where some mathematics gets done. It is said that there's a very short time between being the youngest member of a conference, and becoming the oldest. I don't know quite about Harada-San, but I'm 68, so I'm probably there.

My background is computer science, basically. I started off computing character tables because of a remark by a professor that computing character tables was more of an art than a science, and I thought that that should not be the case, and I was very fortunate in starting off in the 60s, at about the same time as the discovery of the modern sporadics.  Janko's first group was discovered in 1964.

\comment{
I think you should put together a Japanese name for moonshine. I believe you may have done so. It would be a good idea to have a nice national flavor to it. I don't know what the prospects are but 
}
I'm sure you will have some delightful word that would express the contents of moonshine, which means something dubious, among other things, and moonshine, of course, is illicitly produced liquor \footnote{The Kanji, or Chinese characters, for the word is \includegraphics[trim=0mm 0mm 0mm 0mm, clip, width=.5in]{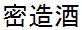}, literally meaning ``secretly made alcohol''. Of course, the word does not quite capture the sense of ``madness'' in English which Conway originally used to express the incredible nature of the Moonshine Conjectures. However, in classical Chinese poetry, numerous allusions are made to drinking accompanied by moon-light. The great poet Li Po (701-762) supposedly drowned himself, in his habitual state of inebriation, trying to grasp the reflections of moonlight in a lake. Thus perhaps \includegraphics[trim=0mm 0mm 0mm 0mm, clip, width=.5in]{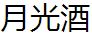}, or ``moon-light liquor'' is a more fitting translation.}.

\subsection{Resources}
\paragraph{Books on Moonshine: }
There is a book by Mark Ronan, the very popular book \cite{ronan}, and there is a series of four fifteen minute talks by him on the BBC3 radio in Europe, at the beginning of this month. They have been recorded. If you want to know more about this recording I can tell you \footnote{
Since its incipience \cite{Con-Nor} and proof \cite{borcherds} (of course, there remains many more things to be understood, including even Ogg's initial mysterious observation of the supersingular primes \cite{Ogg} -- q.v.~recent accounts in \cite{do,sankaran}), Moonshine has developed into a vast field. The reader is referred, for example, to Ronan's book and interview \cite{ronan}, du Sautoy's recent account \cite{Sautoy}, as well as nice technical progress reports of \cite{Gannon:2004xi} and \cite{DGO}.
In parallel, there has been much activity in the physics community extending Moonshine to beyond the Monster, with special focus on Mathieu 24 and its relation to the elliptic genus of K3 surfaces \cite{Eguchi:2010ej,Cheng:2010pq,Cheng:2012tq,Gaberdiel:2010ca,Cheng:2013kpa,He:2012jn,He:2013lha,DGO2}, to matrix models \cite{He:2003pq}, dessins d'enfants \cite{He:2012kw,Tatitscheff:2018aht} and to exceptional Lie algebras \cite{He:2015yoa}.
}.

\paragraph{Using the Web and Moonshine's web page: }
One of the things that I would like to emphasize is the use of the computer and the web today (Internet). I was giving a talk recently and several people come up to me and said: how do you manage to make these astonishing connections between things?
Well, it's not that difficult, really. You have immense resources available on the web, and they grow all the time, and what happened to me, one of the connections I made, was after having got knocked down by a car, maybe again because I was thinking too much about problems. 
I was in bed for two months, and I just searched on the web for phrases in papers, and providing there are not more than a 100 or so papers you can actually go through these papers and see where they are relevant to your interests; a very effective and quite successful way of finding things.

If you really want to know how to get information if you are a graduate student, there is something in the literature: there's Bruce Reznick \cite{reznick} who has written an article on extracting information, and there are some techniques which don't seem to be too well known.

We have had web pages in the past on moonshine groups. One of them was started up by Chi-Han Sah from Stony Brook, and we had quite a nice little group on that, but he died after an operation, and the whole business folded after that. 
Then Chris Cummins, who is a colleague of mine, put up a moonshine page (this was several years ago) with the latest papers and things, and that seems to have disappeared. 
So maybe the time is right to start up something again, where we can discuss things \footnote{
With the growth of blogs, especially in mathematics, there are several emergent sites which are useful \cite{blogs}. However, it would be useful to consolidate these resources, collect comments and have them maintained professionally; much in the spirit of the PolyMath projects \cite{PolyMath}.
}.

\subsection{Talk Outline}
Now, I have three ideas which I think worth pursuing, to the extent that you can show that they are not worth pursuing if necessary, and I will talk about them and explain what little there is to be said about them.

I particularly want to emphasize a few things which are not as well known as they should be. One of them is the action of the Hecke operator and its connection with some very classical objects called \emph{Faber polynomials}. 
Faber was, I believe, a numerical analyst (an analyst), and in 1903, in Mathematische Annalen, he wrote a paper on solving an approximation problem which was of interest to some fairly eminent people, including Hilbert, and that's how
they started. 
But in fact that's not quite true. They go back to a man called Francesco Fa\`{a} di Bruno, and he was an Italian. 
Probably the only Italian mathematical saint. He was beatified in 1988. He died in March 27, 1888.

So these polynomials are of some interest, and one can now look at them from a rather different light. Using these polynomials we can define what we call \emph{replicable functions}. This is a finite class of functions of about several hundred of them,  amongst which the 171 functions which arise in the Monster's context as which we call monstrous moonshine today. These Faber polynomials describe the Hecke action, and that's part of the game.

I spoke to an eminent number theorist a year or two back and he told me that everything was over and we didn't need to think any more about moonshine and we understood everything about it. 
But that's very far from being the case, in my view. 
I think there are several things which are worth thinking about. 
One of them is Witten's idea that there might be some 24-dimensional manifold which would explain this moonshine by looking at the action of the Monster group on the free loop space of the manifold. 
I can't find much by Witten on this, but maybe he's written something. 
That would be a very nice goal to either establish the existence of the manifold and the Monster's action, or to show that such a thing does not exist. 
I think Borcherds, for example, doesn't think that an action on $\mathbb{M}$ exists, but I don't think you should necessarily take much notice of experts; my experience has been rather negative in that respect.

Then finally, as a sort of dream, it would be whether one could gather together all the finite simple groups. 
Let's initially see whether we can put the Monster within a better framework, presumably generalizing the Chevalley work in the Tohoku Journal in 1955 \cite{Chev}, and maybe one can pick up the Monster by generalizing, and it's conceivable that you might be able to pick up other groups; the other six pariahs ($J_1,\ J_3,\ J_4$, Lyons, O'Nan, Rudvalis) in this way. 
John Duncan has found some more moonshine attached to two of the pariahs. 

Let me make another remark. I think that it's quite useful if one finds mathematical objects in other contexts, to find whether there is a connection between them, and I'll say a little bit more about that later. I'll give an example of it shortly.

\subsection{Where to Start ?}
\begin{figure}[!!!t]
\begin{center}
(a)
\includegraphics[width=6in]{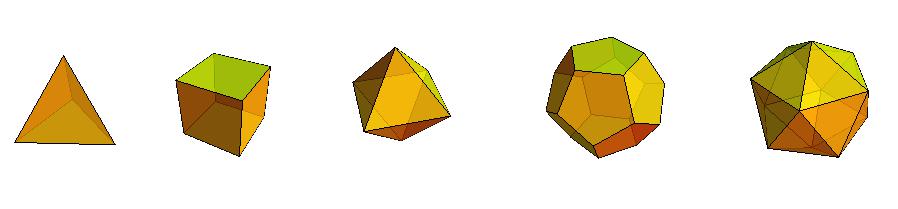}
\comment{
{Table[Graphics3D[{Yellow, Opacity[.8], PolyhedronData[p, "Faces"]}, 
    Boxed -> False], {p, {"Tetrahedron", "Cube", "Octahedron", 
     "Dodecahedron", "Icosahedron"}}]} // TableForm
}
(b)
\includegraphics[width=6in]{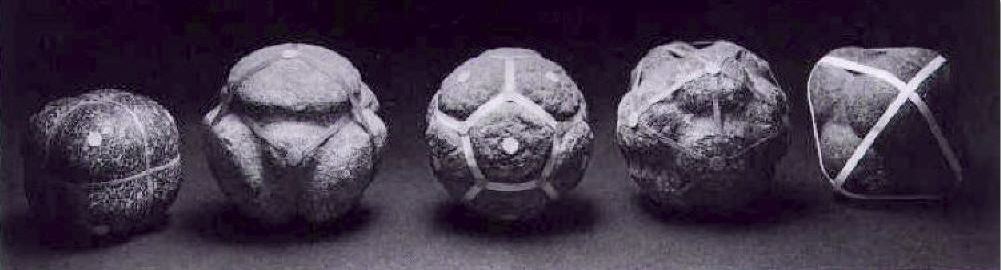}
\caption{{\sf {\small
(a) The 5 Platonic solids, the tetrahedron (T), cubic (C), octahedron (O), dodecahedron (D) and icosahedron (I); C-O and D-I are graph duals and T is self-dual.
(b) The Neolithic carved stones from Skara Brae, Scotland (in the Ashmolean Museum, Oxford), circa 3200BC.
}}}
\label{f:solids}
\end{center}
\end{figure}

\begin{description}
\item[Galois 1832] As of a starting point of the talk, one can start with Galois, who died in 1832, and Galois' work in recognizing the notion of simplicity of a group, normal subgroups, and the other result that $\rm{PSL}_2(p)$ is realizable on the cosets of a subgroup of index $p$ providing that $p$ is not bigger than $11$. 
So there are certain cases of that, which at least initially, were believed to be related to the Monster (see \cite{Con-Nor}).

\item[Mathieu 1861, 1873] One can start with Mathieu in 1861. 
He wrote a paper in 1861 \cite{Mat} in which he said he had found five new groups as transitive extensions of classical linear groups \footnote{
We recall that $k$-transitive means the following. 
Let $G$ be a permutation group on $n$ points and $\{a_1,a_2,\ldots,a_k\}$,
$\{b_1,b_2,\ldots,b_k\}$ are two sets of points with $a_i$ distinct and $b_i$ distinct. If there is an element $g\in G$ mapping each $a_i$ to $b_i$ for $i=1,\ldots,k$, then $G$ is $k$-transitive. 
The only 4-transitive groups are the symmetric group $S_{k \ge 4}$, the alternating group $A_{k \ge 6}$, and the Mathieu groups $M_{24}$, $M_{23}$, $M_{12}$ and $M_{11}$ \cite{cam}.
}.
In the 1861 paper he describes the smaller Mathieu groups $M_{11}$ and $M_{12}$, and he says that others exist. 
Then in a paper 12 years later, he writes that his friends had a bit of trouble seeing how to construct his big groups, and so in 1873 he gives a description of the big groups $M_{24}$, $M_{23}$ and $M_{22}$.

Now, it's notable, and this is true throughout, that the Schur multiplier \footnote{
We recall that for a finite group $G$, the Schur multiplier is the finite Abelian group whose exponent - the LCM of the order of all elements - divides $|G|$.
More generally, the Schur multiplier of a group $G$ is the second group homology $H_2(G; \IZ)$.
} of groups associated with the sporadic groups is larger than one would expect. 
The Schur multiplier of $\rm{PSL}_3(4)$ is an exceptionally big group; a group of order $48$.
Now, $\rm{PSL}_3(4)$ is a group of size $2^6 \cdot 3^2 \cdot 5 \cdot 7 =	20160$ and is the Mathieu group $M_{21}$. 
So $M_{21}$ is {\it not} a sporadic group, but it is a classical group which starts the chain of sporadic groups, $M_{21},\ M_{22},\ M_{23},\ M_{24}$.

I don't know how much skepticism there was when Mathieu wrote this, but there is a paper 25 years after the second Mathieu group by a man called G. A. Miller, who delighted in writing about the problems of other people's work, and his paper attempts to show that $M_{24}$ doesn't exist. 
In 1900 he wrote a paper \cite{Mil} in French correcting himself. So that's a history of the Mathieu groups, and we'll come back to the Mathieu groups later.

\item[Janko 1964] We could start with the Janko groups. 
Janko was a hard worker. He did an enormous amount of work attempting to find sporadic groups, and he
ended up with four groups, which are called $J_1$, $J_2$, $J_3$ and $J_4$. So he worked very hard, and the first successful outcome was in 1964. This is the start of modern era for the sporadic groups.

\item[Plato 400 B.C.] I could start with around the Plato's date, around 400 B.C., with the description of the Platonic solids (cf.~part (a) of Figure \ref{f:solids}). 
Why are we interested in them? 
There's this very curious bijection between the Platonic solids and their symmetry groups inside $\rm{SU}(2)$, and the $A-$, $D-$ and $E-$type Lie structures \cite{ade}. 
So that's the reason for that. 
Predating Plato there's an interesting guy called Empedocles \footnote{
Empedocles (circa 490-430 BC), pre-Socratic Greek philosopher, known as the originator of the cosmogenic theory of the four elements.
}. 
I'll say perhaps a bit more at the end about him. He was interesting because he forecast and predicted the finite speed of light, which I think was quite good for about 500 B.C.

\item[Skara Brae 3000 B.C.] And then there is something less well-known: Skara Brae. 
Skara Brae is a settlement in the Orkneys \footnote{Orkney Islands, northern Scottland, GB.}
which was discovered in about 1850. 
It's called \emph{late neolithic} (that's a cultural date),  but it is about 3200 B.C. 
This settlement, in a very isolated part of the world, contained some carved stones, and I'll show you some pictures of them (cf.~part (b) of Figure \ref{f:solids}). 
Now these stones are all about the same size, and nobody has any idea what they were for.
The belief is they were not weapons because they're not damaged, and the possibility is that they gave the opportunity or permission to speak if you held one of these balls.
I don't know whether the dodecahedron was above the cube or not, but anyhow there are these things around. 
And if you go to the Ashmolean Museum in Oxford \footnote{
Ashmolean collection AN1927.2727-2731, Oxford University,
q.v.~.\url{http://www.ashmolean.org/ash/britarch/highlights/stone-balls.html}
}, they have them there.
What was rather fun was that I mentioned this to Nigel Hitchin, the (emeritus) Savillian Professor of Geometry at Oxford, and is about 400 yards from these things, and he'd never heard of them.
\end{description}

\section{Monstrous Moonshine}
What I would like to do is make some remarks, and see where we get going from here. 
Conway and Norton's paper \cite{Con-Nor} was published at the end of October 1979. 
And the story behind that you've probably all heard, Fischer was visiting me in Montreal, I wrote a letter to Thompson saying that one of the coefficients of the elliptic modular function $j$ was $1$ larger than the dimension of the smallest faithful representation of the Monster. 
Fischer took that back to Princeton. 
I think they all laughed at the concept of there being any connections, but there are. 
Borcherds has reminded me that what had happened was that I was reading a paper by Swinnerton-Dyer and Oliver Atkin \cite{Atk-Swi}, and in that paper they 
give the $q$-expansion for the $j$-function. 

Oliver Atkin used to be a next door neighbour of mine in the ATLAS computer lab, and he is a number theorist.
I was working on these big finite groups, like the Janko group whose order is $175, 560$, 
and he was working on groups associated with two-by-two matrices. I was sure they would be much simpler things than I was working on, but I turned out to be wrong, in retrospect.

This is the order of the monster \footnote{
The Monster, largest of the 26 sporadic finite simple groups, is a 2-generated group, according to the ATLAS \cite{atlas}:
\begin{equation}
M = \langle a,b | a^2 = b^3 = (ab)^{29} = u^{50} = (au^{25})^5 = (ab^2(b^2a)^5b(ab)^5b)^{34} = 1 ; u := (ab)^4(abb)^2 \rangle \ .
\end{equation}
}
\begin{equation}
 |\mathbb{M}|
 =
 2^{46}.3^{20}.5^9.7^6.11^2.13^3.17.19.23.29.31.41.47.59.71  \ .
\end{equation}
There are $15$ primes there. I don't quite know what name they should be given, but anyhow these are the \emph{Monstrous primes}, or the \emph{Monstrous supersingular primes}. 
We will return to these primes shortly.
Thompson makes the remark.
He says the order of a finite group is a very strong invariant. 
These primes appear elsewhere \footnote{
There are marvelous recent expositions on how these 15 primes appear in 5 different contexts by Sankaran \cite{sankaran}, as well as how they can be explained from Moonshine \cite{do}.
} in Erdenberger \cite{Erd} (see later).

That is indeed true, and if you're trying to construct the groups, as we were in the early days, when these sporadic groups were sprouting so that every few weeks there would be a new one, one of the first things was getting hold of the order,  and then using Sylow's Theorem to build up some structure, and perhaps guessing a subgroup, and then using that subgroup and the character table, then using characters building up the character table for the group, and announcing that, and then someone would come along and say the group doesn't exist because the character table doesn't satisfy some property or other. 
Then that property was eventually corrected. 
You had a correct table as far as you knew, and the question was trying to construct the group from the character table, and if that could be done, that was usually done by computer coset enumeration, and then using some technique to prove that the subgroup that you had made exists on the basis of the character table, did indeed exist.

\subsection{Primes in the Monster's Order}
It would be very useful to know more about these primes. 
I don't think there is so much that can be said about them, except for a remark that Ogg made, when he was attending at talk by Serre, I believe, at the Coll\`{e}ge de France \cite{Ogg}. 
This would be in the early 70's. 
One takes
\begin{equation}
\Gamma_0(p)=\left\{
\mat{* & * \\ 0 & *} \bmod \, p 
\right\} \subset PSL(2;\IZ) \ ,
\end{equation}
together with the Fricke involution $\alpha_p:=\tmat{0 & 1 \\ -p & 0}$. 
Consider the group
\begin{equation}
\Gamma_0(p)^+=\left<\G_0(p),\ \alpha_p\right>
\ ,
\end{equation}
and think of it as acting on the upper-half plane \footnote{
Indeed, a classical fact is that the upper-half plane, $\fH := \{z \in \IC : \im(z) > 0 \}$, when adjoining appropriate compactification points known as {\it cusps} which live in $\IQ \cup \infty$, quotients the full modular group to give the Riemann sphere, of genus 0.
}. 
Then the genus of the Riemann surface $\G_0(p)^+ \backslash \mathbb{H}$ is zero precisely when $p$ is one of the 15 supersingular primes that appear in the Monster's order. 
That's one number theoretic characterization of these primes \footnote{
Genus 0 congruence subgroups are very rare.
For example, there are only 33 which are torsion free \cite{classSebbar} and the relation of these with elliptic surfaces, especially with K3 surfaces, is discussed in \cite{sebbar,mckaysebbar,He:2012kw,He:2012jn}.
}. 
We still don't know why, and Thompson regarded that as one of the major
questions to be answered in connection with the Monster.

There is another way of saying it: for elliptic curves defined in characteristic $p$, then all the supersingular $j$-invariants of these curves (being a priori in $\mathbb{F}_{p^2}$) are lying in the base field $\mathbb{F}_p$, rather than in $\mathbb{F}_{p^2}$ , precisely if $p$ is one of the above 15 primes.

Now, rummaging through the contents of the preprints on the arXiv.org every weekday you look through and see if there's anything of interest. 
We found a paper by Cord Erdenberger \cite{Erd}, who is a student of Klaus Hulek from Hannover. 
And these 15 primes come up in his work. 
His title is ``The Kodaira Dimension of Certain Moduli Spaces of Abelian Surfaces'' (MR20923323 (2004)). He considers Abelian surfaces $(1,p)$ polarized ($p$ a prime), and uses Jacobi cusp forms of weight $2$ and level $p$, and these apparently exist just when $p$ does not divide the order of the Monster. 
So in a sense, one might say that they are related to 
the existence of this Monster group. 
This is one connection which needs some explanation.
Here you are working with a subgroup of the symplectic group rather than the modular group. 

I've contacted Erdenberger and his supervisor Hulek, and nobody seems to know quite whether this is really saying something new, or whether it can be interpreted in terms of these supersingular elliptic curves that I mentioned earlier. It's something that should be pursued, at least try to find out whether there is a connection or not.

I would suggest if you want to follow this up, look in Math review for the Math Review number I have given above. 
The reason being that Sankaran reviewed it, and he does mention this in the review. 
You won't find any paper about it, and certainly Erdenberger was not aware of it. 
I don't know whether anyone is pursuing this; 
I don't know of any pursuit of this fact \footnote{
The reader is referred to the recent works of \cite{do} and  \cite{sankaran} for various explanations.
}.

\subsection{Balance}
This is typical of the sort of fact that you can gather, and one can formalize it as something (I don't know if it's a great thing to do so) maybe the words are ``retro-syntactic retrieval'', or something like that, but the game is very simple \footnote{
This was mentioned earlier in the introduction about how one could retrieve information and establish correspondences.
}. 
You have a bunch of people working on different subjects, and then if you study the phrases that are used in common by these people, or you find people who use a common phrase, there's a good chance that there is some related activity going on between the people that use this phrase. 
I didn't do that in this case. I think I was just looking through the arXiv.org and found it.

Now let me say something about balance. If we go back to the first paper, the word used is ``seminal'', certainly it was the only paper for a long time on the subject by Conway and Norton \cite{Con-Nor}. 
Conway is here, and Simon Norton is not here. 
I don't think he should be forgotten. 
He  is very often the motivating force between a lot of activity, some of which never gets published.

In this paper there is a list of observations which are introductory to the business of moonshine and one of them is that elements of the group $M_{24}$ are \emph{balanced}. 
So, what does this mean, and what is its significance? 
If you take a permutation in terms of disjoint cycle lengths you have a bunch of numbers which form a partition of the degree. 
A permutation is \emph{balanced} if the product of the lengths of pairs from the outside-in is constant. 
Here is an example \footnote{
Thus, an element of the permutation group $S_{24}$ of degree 24, would have cycle notation $(a_1)(a_2, a_3)(a_4, \ldots, a_{10})(a_{11}, \ldots, a_{24})$, which is indeed the shape of one of the $1575$ conjugacy classes of $S_{24}$ (1575 is the number of unrestricted partitions of 24).
}, a permutation with cycle lengths $1$, $2$, $7$ and $14$ is balanced, since $N=1.14=2.7=14$. 
The number $N$ is called the \emph{balance number}, if it exists.

In 1980 or thereabouts, we had a conference called ``The Coming of Age of the Finite Groups'', and some of the people were here then. Dummit and Kisilevsky and myself classified  all permutations of degree $24$ that are balanced \cite{DKM}. 
Why choose 24? Well, we're going to replace $k$ by $\eta(q^k),\ q=e^{2\pi i\tau}$, for each cycle of length $k$, and thus form the product, we call it \emph{eta products}.

We found all the $\eta-$products which are weakly multiplicative in the coefficients\footnote{
The simplest case is the famous $\Delta(q) = \eta(q)^{24}$, which is the modular discriminant function, with q-expansion $\Delta(q)=\sum\limits_{n =1}^\infty \tau(n) q^{n}$ with $\tau(n)$ being the Ramanujan tau-function.
This is weakly multiplicative in the sense that
$\tau(m \ n) = \tau(m) \tau(n)$ if $\gcd(m,n) = 1$.
The multiplicative eta-products appear in physics, especially in partition functions in string theory and are discussed in \cite{Govindarajan:2009qt,Cheng:2010pq,He:2013lha}.
}. 
There are exactly $30$ of them, and all the permutations in the Mathieu group $M_{24}$ are balanced. 
And being balanced and weakly multiplicative is of the same thing as a theorem of Bryan Birch and Morris Newman on that \cite{New}.
And you have a cusp form for each balanced permutation, with what's called a
{\em Grossen-character}, whose weight is half the number of parts.

Now we can generalize this to eta-quotients instead of eta-products by writing fractions. 
Here is an example: $2^{24}/1^{24}$ which means  $\eta(2\tau)^{24}/\eta(\tau)^{24}$.
The multiplicative $\eta$-quotients have been classified too \footnote{
Cf.~also \cite{kilford,MarOno} for relations to elliptic curves.
} by Yves Martin \cite{Mar}.
These appear in \cite{Apo}. 

The products are straight-forward because there are only finitely many partitions to look at, so you just go through them, find out the ones which look like they're multiplicative by looking at the first few coefficients and checking, and then filter them out and then you have to prove something.
But for the quotients that's a different matter. 
The quotients are much more difficult. 
There is potentially infinitely many of them. 
The guy who had done a paper on them is Yves Martin. 
He hasn't completely done it. 
He made an assumption that both the eta-quotient and the eta-quotient with the above Fricke involution action on it are weakly multiplicative and that's not asked for. 
So, the general question about what eta-quotients are multiplicative is not known.

Being multiplicative means that you have some Euler product through the (inverse) Mellin transform \footnote{
That is, we can form the Dirichlet series for the coefficients $c_n$ of the q-expansion $\sum\limits_n c_n q^n$ for these eta-products to give $\eta_g(s) = \sum\limits_n c_n n^{-s}$.
This can then be taken as product over primes as in \eqref{euler}.
}, and here is the Euler product (we use the subscript $g$ to identify the particular eta-product):
\begin{equation}\label{euler}
\eta_g(s)= \prod_p \left(
     1-\frac{a_p(g)}{p^s}+\frac{b_p(g)}{p^{2s}}
     \right)^{-1}
\end{equation}
with $a_p$ and $b_p$ integers; we have
\begin{equation}
a_p^2-a_{p^2}=b_p=\chi\left(\frac{\cdot}{p}\right)p^{w_g-1} \ ,
\end{equation}
where $\chi$ is some character and $w_g$ is the weight of the eta-quotient, which is half the number of parts. 
In particular, on the identity, the partition is $1^{24}$ and the weight is $w_g=24/2=12$, which yields $p^{11}$. 
We have a degree $p^{11}$ generalized character. 
These are proper characters which can be checked this directly. There  is some anomalous behavior for $p=3$, but other than that everything is clear.

\subsection{Conjugacy Classes}
There are curious connections with physics and conjugacy classes, 
and I think someone raised this the other day at the EWM meeting \footnote{
Encounters with Mathematics, Chuo University, May, 2008,
\url{http://www.math.chuo-u.ac.jp/ENCwMATH/45.shtml}.
}, 
that in studying elliptic genera you look at commuting pairs of elements, and that's related to the number of conjugacy classes in the group when working with a finite group. 
And I just make the passing remark, 
the class number - the number of conjugacy classes- of $M_{24}$ is 26, which should ring a bell with some physicists \footnote{
In string theory, the critical space-time dimension of the bosonic string is 26 and that of type II super-string is 10.
}. 
You have 5 quadratic boxes of irrationalities, so there are 21 classes of cyclic subgroups, or if you like, 21 rationally irreducible representations. 
The group $M_{12}$, believe it or not, has $10$ instead of $26$.

If you want to know more about this remark about conjugacy classes, there is recent publication, a large book in fact called "From Number Theory to Physics", in Les Houches proceedings of 2002, containing a paper with Sebbar and myself
\cite{McK}. 
More recently, in a paper in a conference proceedings on things to do with moonshine, there's a paper by Anda Degeratu and Katrin Wendland \cite{Wend},
and they look again at one of these situations that you discover by reading and saying ``my goodness, it's the same number there.'' 
They are looking at a situation where the appropriate number replacing 26 is 194, and 194 is the number of conjugacy classes in the Monster \footnote{
In other words, as remarked in \cite{McK}, the number of conjugacy classes of $M_{12}$ and $M_{24}$ are respectively 10 and 26, the critical dimension of the supersymmetric and the bosonic string theories.
Moreover, 194, the number of conjugacy classes of the Monster, is the Picard number of the base of an elliptically fibred Calabi-Yau threefold in an extremal case of heterotic-F-theory duality as studied in \cite{Aspinwall:2000kf}.
Furthermore, of these 194 classes, considered as column-vectors in the character table, only 163 are linearly independent; and of course, 163 is a famous Heegner number where the exponential assumes an almost-integer value: $\exp(\pi \sqrt{163}) \sim 640320^3 + 744$.
}. 
So they are looking at a situation that is  interesting if anything comes out of it.

\subsection{Frame Shape}
Now, these shapes we've been talking about -- the partition of $n$ where $n$ is the degree, or the more general situation when you divide one term by another is called \emph{Frame shapes}. 
Now ``Frame'' is the name of a person, J.~S.~Frame, not an abstract notion, and he was an interesting man. 
His thesis was on character tables in the 1930's \cite{frame}. 
He did character tables like other people do crosswords. 
So that's who Frame was. 

Basically speaking, what you are doing is that you are describing the eigenvalues within an orthogonal group. 
So if you have a fraction, you want to make sure that the eigenvalues that you take away from the denominator are already in the numerator, to make any sense. 
There's a paper by Takeshi Kondo \cite{Kon}, who wrote a very nice paper about the Frame shapes of elements in the automorphism group of the Leech lattice, Aut$(L)$.
If you go down from there to Conway's group $\cdot 1$, which is this group modulo the action on diameters - i.e., Aut$(L) / \left< \pm1 \right>$ - you get a mixture of the functions that describe the elements on the various classes and they are not as consistent as those for the monster $\mathbb{M}$.

By the way, there's a very nice survey by Masao Koike \cite{koike}, I will say more about him later on.
Again this was in Sugaku, in Japanese that has been translated into English by the AMS translations $\#160$.
This is one of the few early surveys on moonshine, so that's a useful paper too.

\subsection{Faber Polynomials}
If we take a Riemann map \footnote{
This approach of looking at Moonshine from the perspective of geometric function theory, in terms of the shape of the analytic functional form of $j(q)$ and generalizations, is very much the spirit of the current lecture notes, and is also summarized in \cite{essentials}.
} from the exterior of some region in the complex plane containing two points at least, by the Riemann mapping theorem we can map the exterior of this region with some conditions at infinity to the exterior of a disk of radius $d$, and we can normalize this
\begin{equation}\label{riemannmap}
z=\phi(w)=dw+d_0+\sum_{k\geq 1}d_kw^{-k} \ .
\end{equation}
If you are doing analysis you don't have to worry about the constant term and if you are doing moonshine you put it equal to 0 and the radius equal to 1.
We get an inverse of the same forms as \eqref{riemannmap}
\begin{equation}
w=\phi^{-1}(z)=z/d+g_0+\sum_{k\geq 1}g_kz^{-k},
\end{equation}
and similarly we can take $d=1$ and $g_0=0$.
The {\em Faber polynomial} $F_n$ is the part of $(\phi^{-1}(z))^n$ with non-negative powers of $z$. So you're picking up the polynomial part of this series $(\phi^{-1}(z))^n$.

You may not be familiar with this, but if you're looking at pseudo-differential operators that's a standard procedure to pick the plus part of the operator. That's what the Faber polynomial does for you, and that how it's used and I'll say quite a lot more about that \footnote{
In other words, we consider a meromorphic function and its inverse with Laurent expansion of the form \eqref{riemannmap} and \eqref{eqn:f}.
This is clearly inspired by the form of the q-expansion of the $j$-invariant, as we shall shortly see.
We emphasize that the constant term is 0, so henceforth, by the $j$-invariant, we mean the normalized one $j(q) - 744$.
In \cite{He:2003pq}, this shape of a Laurent series was interpreted as the master-field of a large N matrix model, whereby giving a modular matrix model.
}.

What we do is we compose with the map $z\mapsto 1/q=e^{-2\pi i z}$, to get
\begin{equation}\label{eqn:f}
f(z)= q^{-1}+\sum_{k\geq 1}a_kq^k,\ \quad \Im(z)>0\ .
\end{equation}
The $a_k$'s are general coefficients and we'll take them to be integers, but they need not be integers, and Simon Norton has classified the functions we are interested in (the replicable ones to come later) even when they have complex coefficients. 
If these coefficients are not integers, they do lie in a field whose Galois group is an elementary $2$-group over the rationals. 
In other words, the $a_k$'s, lie in a composite of quadratic fields.
I don't think we really know which quadratic fields and why, but anyhow that's where they lie when we're talking about replicable functions, that's the ones we are interested in.

But for our purposes and for all the stuff here we work with the $a_k$'s being integers. That avoids any problems with Galois theory and is convenient. The functions of the form \eqref{eqn:f} are typical functions we shall study, and we are going to study them first of all slightly more generally than the connection with the Monster, and then specialize to functions that are attached to $\IM$.

\subsection{Grunsky Coefficients}
We can define the elliptic modular function $j$ by the above property, because, if this holds for some other function $f(z)$ for all positive integers $n$, the level of this function $f(z)$ must be equal to $1$, and is therefore a rational function of $j(z)$, and so, using the fact that $j$ is normalized at infinity as in \eqref{eqn:f}, we have that $f(z)=j(z)$.

For $f$ as in \eqref{eqn:f} we write
\begin{equation}\label{hnm}
F_n(f)=q^{-n}+n\sum_{m\geq 1} h_{m,n}(f)q^m \ .
\end{equation}
the coefficients $h_{m,n}$ are called \emph{Grunsky coefficients}, see \cite{Grun,Pomm}. These are in fact symmetric in the $m,n$ indices
\footnote{\label{ft:Fn}
We can in fact define the Grunsky coefficients and Faber polynomials in the following way. We will encounter some of the ensuing expressions in due course.
Let $g(z)$ be a holomorphic univalent (i.e., one-to-one on the open set) function on the unit disk $|z| < 1$, normalized so that $g(0)=0, \ g'(0) =1$.
Then the function $f(z) = g(1/z)^{-1}$ is a non-vanishing univalent function outside the unit disk with simple pole at $\infty$ with residue 1; that is,
\[
f(z) = z + a_0 + a_1 z^{-1} + a_2 z^{-2} + \ldots \ .
\]
The expansion coefficients $c_{nm}$ of
\begin{equation}\label{cnm}
\log \frac{f(z) - f(w)}{z-w} = - \sum\limits_{m,n>0} c_{nm} z^{-m} w^{-n}
\end{equation}
are the {\em Grunsky coefficients}.
Definition \eqref{cnm} implies, upon $z \frac{\partial}{\partial z}$ on both sides, that $\frac{zg'(z)}{f(z) - f(w) } - \frac{z}{z-w} = \sum\limits_{m,n>0} m c_{nm} z^{-m} w^{-n}$. Thus we define the {\em Faber polynomials} $F_n(w)$, as
\begin{equation}\label{Fn}
\frac{zg'(z)}{f(z) - w } = \sum\limits_{n \ge 0} F_n(w) z^{-n} \ .
\end{equation}
It is non-trivial that, thus defined, $F_n(w)$ are monic polynomials of degree $n$.
In fact, $F_n$ are themselves polynomials in the coefficients $a_i$ in the definition of $f(z)$.
This is seen as follows.
Definition \eqref{Fn} implies, upon applying $\int_0^\infty dz \frac{1}{z}$ on both sides, that $\log \frac{f(z) - w}{z} = - \sum\limits_{n \ge 1} \frac{1}{n} F_n(w) z^{-n}$.
Expanding out $f(z)$ order by order and comparing with \eqref{Fn} then gives the recursion
\[
F_n(0) = 1 \ , \quad
F_n(w) = (w-a_0) F_{n-1}(w) - n a_n - \sum\limits_{i=0}^{n-1} a_{n-i} F_i(w) \ .
\]
Furthermore, combining \eqref{Fn} and \eqref{cnm} we have that
$\sum\limits_{n \ge 0} F_n(g(\zeta)) z^{-n} = \frac{z}{z-\zeta} + \sum\limits_{m,n >0} m c_{nm} z^{-m} \zeta^{-n} = \sum\limits_{n>0} \left(\frac{w}{z}\right)^n +
\sum\limits_{m,n >0} m c_{nm} z^{-m} \zeta^{-n}$ so that
\[
F_n(g(z)) = z^n + \sum\limits_{m \ge 1} c_{nm} z^{-m} \ ,
\]
which is \eqref{hnm} in our definition, up to the factor of $n$ which will be more convenient for our succeeding discussions.
}.
They have a remarkable connection with the Bieberbach conjecture
\footnote{
We recall the statement of the Bierberbach Conjecture, proven by de Branges.
For univalent holomorphic function with Taylor series of the form
$f(z) = z + \sum\limits_{n\ge2} a_n z^n$ (such functions are called {\it Schlicht}, or simple/plain), the coefficients have the property that
$|a_n| \le n$ for all $n \ge 2$.
}.
The Bieberbach conjecture is a bound on the coefficients of functions univalent on the unit disk, and there's a very nice book \cite{de bran}, an AMS publication, on the solution to the Bieberbach conjecture by de Branges, and these Grunsky coefficients played the major role in the establishment of this conjecture in its early days\footnote{
There is a very recent paper on the appearance and relevance of Bieberbach/de Brange as well as Grunsky coefficients in scattering amplitudes in quantum field theories \cite{Haldar:2021rri}.
}.

In our definition \eqref{hnm}, we have written $F_n(f)$ in this way with an $n$ in front in order to take advantage of the symmetry of $h_{m,n}$ in $m$ and $n$. In fact that is the definition of the Grunsky coefficients
\begin{equation}
h_{m,n}(f)=T_n(f)|_{q^m} \ .
\end{equation}
There will be other introductions to the Faber polynomials and these Hecke operators $T_n$ later \footnote{\label{ft:Faber}
In \cite{McK}, a particularly nice characterization of the Faber polynomial is as given.
Consider, as always, a function $f(q) = q^{-1} + \sum\limits_{n\ge1} a_n q^n$ for nome $q = \exp(2 \pi i z)$ with $\im(z)>0$, as in \eqref{eqn:f}.
Then, for each $n \in \IZ_{>0}$, there is a {\it unique} monic polynomial $F_n$ such that 
\[
F_n(f(q)) = q^{-n} + \cO(q) \ , \quad \mbox{ as } q \rightarrow 0 \ . 
\]
These are the Faber polynomials.
Depending on the Taylor series of $f(q)$, the first few are
\[
F_0(z) = 1, \
F_1(z) = z, \
F_2(z) = z^2 - 2 a_1, \
F_3(z) = z^3 - 3a_1 z - 3 a_2 \ .
\]
More generally, we have
\[
F_n(z) = \det(z \II - A_n) \ , \qquad
A_n := \tmat{
a_0 & 1 &&&&
\\
2a_1 & a_0 & 1 &&&
\\
\vdots & \vdots & \vdots &&
\\
(n-2)a_{n-3} & a_{n-4} & a_{n-5} & \ldots & 1 &
\\
(n-1)a_{n-2} & a_{n-3} & a_{n-4} & \ldots & a_1 & 1
\\
na_{n-1} & a_{n-2} & a_{n-3} & \ldots & a_1 & a_0
}
\]
}.

It's not obvious from the above expression that $h_{m,n}$ is symmetric in $m$ and $n$, but it is, and we can see it with a slightly different generating function for the Grunsky coefficients.


\section{Hecke operators and Faber Polynomials}
I now turn to Part 2 of my lectures, having alluded to the Hecke operators.
There are Hecke operators, and often in the books they assume that the function on which the Hecke operator is acting has a weight greater than zero, whereas the functions we are interested in all have weight zero. 
You have a group action acting linear fractionally on $z$, and on the modular function $j$ this is given by,
for all $n\geq 1$
\begin{equation}\label{eqn:Tn}
     T_n(j(z))
     =
     \frac{1}{n}
     \sum_{\substack{ad=n\\0\leq b<d}}
     j\left(
     \frac{az+b}{d}
     \right)
     =\frac{1}{n}F_n(j(z))\ .
\end{equation}
The effect of the Hecke action is to replace the pole of order 1 of the $j$-function with a pole of order $n$ at infinity. 
Moreover, the action of the Hecke operator preserves the space of modular functions. 
Hence, $T_n(j(z))$ is a rational function of $j(z)$, and so this going to be a polynomial in $j(z)$. 
$T_n(j(z))$ can be expressed both as a $q$-series and as a polynomial in $j(z)$. In fact, we can define the $j$-function by the fact that there is an action of the Hecke operator defined in terms of sum over the function valued on sublattices. Let's have a quick look at this, which is standard.

\begin{figure}[!h!t!b]
\begin{center}
(a)
\includegraphics[width=1.5in]{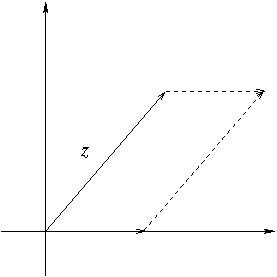}
\qquad $\stackrel{T_2}{\longrightarrow}$ \qquad 
(b)
\includegraphics[width=3in]{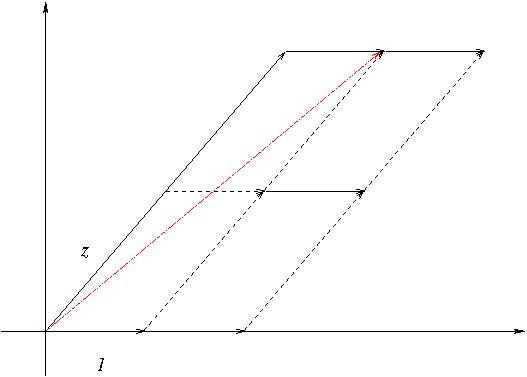}
\caption{{\sf {\small
(a) Fundamental region of the lattice $\Lambda = \IZ + z\IZ$.
(b) The 3 sublattices of index $n=2$, viz., $\IZ + 2z\IZ$, $\IZ + (1+2z) \IZ$ and $2\IZ + z \IZ$.
}}}
\label{f:lattice}
\end{center}
\end{figure}

In part (a) of Figure \ref{f:lattice}, this is the fundamental region, up to orientation and homothety.
The Hecke operator $T_n$ maps the lattice $\Lambda=\BZ+z.\BZ,\ z=\omega_1/\omega_2$ such that $\Im z>0$, to sublattices $\{\Lambda_i\}$ of index $n$, so induces an action on functions defined on these lattices 
\footnote{
The fundamental fact here is that sublattices of index $n$ are in one-one correspondence with integer matrices $\tmat{a & b \\ 0 & d}$ with $a > 0$, $b = 0, 1, \ldots, d-1$ and $a d = n$.
For example, at $n=2$, we have 3 such matrices,
$\tmat{1&0\\0&2},\tmat{1&1\\0&2}$ and $\tmat{2&0\\0&1}$, corresponding to the lattices $\IZ + 2z\IZ$, $\IZ + (1+2z) \IZ$ and $2\IZ + z \IZ$.
}. 
For example, for $n=2$, there are 3 lattices of fundamental region twice that of $\Lambda$ and the Hecke action is given by
\begin{equation}
T_2:\ f(z)\mapsto\ \frac{1}{2}\left(f(2z)+f\left(\frac{z}{2}\right)+f\left(\frac{z+1}{2}\right)\right) \ .
\end{equation}
This is drawn in part (b) of Figure \ref{f:lattice}.
I've used $z$ here because $\tau$ is used later.

For each of the functions we are interested in there will be a discrete subgroup $G_f$ of $\PR$ with respect to which $f$ is modular, and in \cite{Con-Nor} there is a discussion on how much is being fixed. The action of $G_f$ is linear fractional.
For us, the transformation as an element of $G_f$
\begin{equation}\label{Gf}
f\left(\frac{az+b}{cz+d}\right)=(cz+d)^kf(z)
\end{equation}
is not affected by the automorphy factor $(cz+d)^k$ because we are working with functions of {\em weight} $k=0$, and this is rather important since the behavior is different \footnote{
That is, the $j$-function is an absolute invariant
$j\left(\frac{az+b}{cz+d}\right)=j(z)$ for $\tmat{a&b\\c&d} \in PSL(2;\IZ)$.
Indeed, for weight $k$ objects, the standard definition \cite{Serre} of the Hecke operator is
\[
T_n(f(z))
     =
     n^{k-1}
     \sum_{\substack{ad=n\\0\leq b<d}}
     d^{-k}
     f\left(
     \frac{az+b}{d}
     \right)
\ .
\]
for all $n \in \IZ_{\ge 1}$.
}.
There is also the notion of the {\em level} , which is the smallest $N$ such that $\G(N)\subset G_f$ 
\footnote{
The congruence groups are defined with some modulo $N$ relation.
For example, the principal congruence subgroup is 
$\Gamma(N) := \{A \in PSL(2;\IZ)  \, \left| \, A \equiv I \bmod N \right.\}$.
}.

Classically, for all $\gcd(n,N)=1$, $T_n(f)$ is a polynomial in $f$ ($f$ as in \eqref{eqn:f}). 
This last statement is the critical one that enable us to generalize the action of the Hecke operator, and the generality is that if one defines the $j$-function with the normal Hecke operator acting as a polynomial then what we do is we preserve the action of the (Faber) polynomial and redefine the Hecke operator. Indeed, since for all $n\geq1$ we have
\begin{equation}\label{eqn:Tnf}
T_n(f(z))=\frac{1}{n}
     \sum_{\substack{ad=n\\0\leq b<d}}
     f\left(
     \frac{az+b}{d}
     \right)
     =\frac{1}{n}F_n(f(z))\ ,
\end{equation}
the modular level is one, and therefore is a rational function of $j$. 
Hence, using \eqref{eqn:f} where $f$ is normalized to have a simple pole at infinity and a zero constant term, implies that $T_n(f(z))$ is a polynomial in $j(z)$, and thus $f=j$.
This enables us to define the $j$-function from the property that the Hecke operator on $j$ averaging over sublattices of index $n$ is the Faber polynomial associated with $j$. 

I want to say quite a bit about this important polynomial. 
It's easy to find this polynomial from a computational and algorithmic points of view. 
Start off with \eqref{eqn:f} and look at various powers of $f$ you'll get various negative powers of $q$ on the right hand side of \eqref{eqn:Tnf}, and by forming linear combinations of these powers of $f$ on the left hand side of \eqref{eqn:Tnf} we can eliminate all but the largest negative power of $q$, and if there is a constant term we can put it into the polynomial on the left hand side. 
That makes it clear that this polynomial is very simply determined from \eqref{eqn:f} and is unique. 
What this polynomial is doing is replacing the simple pole at infinity in \eqref{eqn:f} by an order $n$ pole also at infinity in \eqref{eqn:Tnf}. 
And if we write, for $f$ as in \eqref{eqn:f},
\begin{equation}\label{eqn:h_m,n}    
F_n(f)=q^{-n}+n\sum_{m\geq 1} h_{m,n}(f)q^m.
\end{equation}
the coefficients $h_{m,n}$ are the \emph{Grunsky coefficients}. 
They have a lot of interesting properties too
\footnote{
Let us re-iterate this point.
We saw in the footnotes above that the Faber polynomials are the unique degree $n$ monic polynomials bringing $q^{-1} + \cO(q)$ to $q^{-n} + \cO(q)$ as $q \rightarrow 0$.
Now, our Hecke operator on $j$ of weight 0 as in \eqref{eqn:Tnf}, thus all $T_n(j(z))$ are invariant under $PSL(2;\IZ)$ since the sum gets permutated by the action of the modular group. Whence they must be rational functions in $j(z)$ since the $j$-invariant, being a Hauptmodul \footnote{
I want to use the word, rather than call things Hauptmodul because Hauptmodul has a nasty plural, and we would like to try to call them {\em principal modulus}, and we have done that in some papers.
}
, generates the function field of invariants. However, since it has no poles in the upper half plane, they must in fact be polynomials.
In fact, we find that $T(j(q)) = q^{-n} + \cO(q)$ as $q \rightarrow 0$.
By uniqueness then, $T_n(j)$ must be (up to overall normalization) the Faber polynomials!
}
.

\subsection{Replicable Functions: Norton's Basis}
We can define replicable functions $\{f^{(m)}\}$, say, by the same $q$-expansion as in \eqref{eqn:f}
\begin{equation}\label{eqn:f^(m)}
f^{(m)}(z)=\frac{1}{q}+\sum_{k\geq1}h_k^{(m)}q^k\ ,
\end{equation}
I am using the superscript $(m)$ here, and I was talking to John Conway in the breakfast about this. 
I think the notation has to be changed. I've been using small letters for coefficients and character, and I've used a little $h_k$ here instead of $a_k$. You can use $a_k$ if you wanted to. But the I use the $h_k$ to remind you of the characters of the Monster in the special case when you restrict to the Monster. We can define a collection of functions here by, for all $n\geq1$
\begin{equation}\label{eqn:T^nf}
\widehat{T}_n(f)=\frac{1}{n}\sum_{\stackrel{ad=n}{0\leq b<d}}f^{(a)}\left(\frac{az+b}{d}\right)=\frac{1}{n}F_n(f)\ .
\end{equation}
What we've done is that we kept  the Faber action in the above and replaced the sum over representatives of the sublattices of the function $f$ by the functions $f^{(a)}$. Those of you who want a glimpse of the future, I can give you the relation between $f$ and $f^{(a)}$ with reference to the Monster. 

To an element of the Monster $g$ there corresponds a function $f_g$, corresponding for the moonshine for this function on the element $g$, then raising $f$ to the $a-$th replicate power $f^{(a)}$ is the same as replacing $g$ by $g^a$. You will get all this in a short time. 
There is no need to make reference to the Monster in order to define this \footnote{
Historically, the concept of replicability came about from Conway-Norton's initial observation \cite{Con-Nor} that the moonshine functions (McKay-Thompson series) obeyed certain functional identities, which they called replication.
This is the reason for studying the type of recurrences in \S\ref{s:mahler}.

The motivation in defining it in the manner of the present section is to generalize the remarkable fact that action of the $n$-th Hecke operator on $j$ is the $n$-th Faber polynomial in $j$.
Thus a function of the expansion type \eqref{eqn:f} is {\em replicable} if there exists a family of function $\{ f^{(a)} \}$, called replicable functions of $f$ such that the generalized Hecke operator on these $\frac{1}{n}\sum\limits_{\stackrel{ad=n}{0\leq b<d}}f^{(a)}\left(\frac{az+b}{d}\right)$ is the Faber polynomial in $f$, i.e., $\frac{1}{n}F_n(f)$.
}.

There is an inductive definition, and there is only going to be one term, which is $f^{(n)}$, in the sum in \eqref{eqn:T^nf}, and so you can take out the rest of this sum and put it on the right hand side and that with will involve only $f^{(a)}$ with $a<n$ and you have an inductive definition of these functions. I'll say more about these Hecke operators.

Now, Norton did all this unaware that it had all been done 70 years before and before that.
Perhaps the easy definition of the Grunsky coefficients $h_{m,n}$ is in terms of this generating function
\begin{equation}\label{eqn:h_m,n-Norton}
\ln\left(\frac{f(p)-f(q)}{\frac{1}{p}-\frac{1}{q}} \right)=-\sum_{m,n\geq1}h_{m,n}p^nq^m
\end{equation}
with $p=e^{2\pi is},\ q=e^{2\pi iz}$.
You take the logarithm of the difference of the function evaluated at two different arguments $p$ and $q$. Remember that $f(p)$ (resp. $f(q)$) will start with $1/p$ (resp. $1/q$), so dividing through by $1/p-1/q$ gets rid of the singularities and you end up with a quite nice series, namely the series on the right hand side of \eqref{eqn:h_m,n-Norton}. That's the generating function for the $h_{m,n}$ and you can expand that to have
\begin{equation}\label{eqn:ln-expqnsion}
Eq.~\eqref{eqn:h_m,n-Norton}=
\ln\left(1-pq\sum_{k\geq1}h_k\frac{p^k-q^k}{p-q}\right)
\end{equation} 
and there are certain consequences of the expansion one of them is that the $h_{m,n}$'s are polynomial in the $h_k$'s with $k\leq m+n$; the value $m+n$ provides you with the grading. I'm using the convention $h_k=a_k=h_{k,1}$. When we compare coefficients in the last expansion with the $\ln$ term, we get this recursive expression
\begin{equation}\label{eqn:h_m,m-recursive}
h_{r,s}=h_{r+s-1}+\frac{1}{r+s}\sum_{m=1}^{r-1}\sum_{n=1}^{s-1}h_{m+n-1}
(r+s-m-n)h_{r-m,s-n}.
\end{equation}

Indeed, in \eqref{eqn:ln-expqnsion}, the terms $\displaystyle\frac{p^k-q^k}{p-q}$ are part of $p^sq^t$ where $s+t=k-1$, but you are multiplying by $pq$ so you get a term involving the expression involving the $h_k$'s and the $h_{m,n}$ corresponding to the appropriate exponent. If we call $r+s$ the grade then all $h_{r-m,s-n}$ in \eqref{eqn:h_m,m-recursive} have a lower grade. 
The $h_{r,s}$ are not integers, but the largest denominator is $(r,s)$. 
I believe that's correct in a sense. 
Suppose we want to compute a $q$-coefficient of the function, then that is typically given by the term $h_{r+s-1}$. So we have a choice: if we want to compute $h_{k}$ then we can choose $r+s=k+1$ so that $k=r+s-1$, and we can do that what will give us the initial step, and the game is to try to find a pair which has a reduced sum. 

Simon Norton \footnote{
Following \cite{McK}, we can proceed with this formal definition of a replicable function.
Consider a function the form \eqref{eqn:f}, and write its corresponding Faber polynomial, with Grunsky coefficients $h_{m,n}$ as in \eqref{hnm}.
Then $f$ is replicable if $h_{m,n}=h_{r,s}$ whenever $\gcd(m,n) = \gcd(r,s)$ and ${\rm lcm}(m,n) = {\rm lcm}(r,s)$.

Equivalently, we can define replicable functions using the Hecke operators (for weight 0).
The function of our form \eqref{eqn:f} is replicable, if for each positive interger $n$ and positive divisor $a | n$, there are functions $f^{(a)}$ of the form of \eqref{eqn:f} such that
$\widehat{T}_n(f) := F_n(f(q)) = \sum_{\stackrel{ad=n}{0\leq b<d}}f^{(a)}\left(\frac{az+b}{d}\right)$.
The functions $f^{(a)}$ are called replication powers and have the property that
\begin{equation}
f^{(k)}(q) = q^{-1} + \sum\limits_{i \ge 1} \left( 
k \sum\limits_{d | k} \mu(d) h_{dki, \frac{k}{d}}
\right) q^i \ ,
\end{equation}
where $\mu$ is the standard M\"obius $\mu$-function.
}
has a definition of a replicable function which is that a function is replicable if $h_{m,n}=h_{m',n'}$ whenever $mn=m'n'$ and $\gcd(r,s)=\gcd(m,n)$. This is an important definition and we can take advantage of this in computing the coefficients of a replicable function from \eqref{eqn:h_m,m-recursive}. 

The following picture shows how to compute the coefficient of $q^k$. 
With $k$ fixed this gives us a choice of $r,\ s$, so we can draw the line $x+y=r+s=k+1$, and then the game is to find a point $(r,s)$ on that line that dominates, if we are lucky, some other point $(r',s')$ with the same $\gcd$ and $\rm{lcm}$. We then have the hyperbola $xy=rs$.
So if we can find this $(r',s')$ on a lower line then we start with, we can go to this line and proceed to do the same thing again, and each time we do that there are two possibilities --  there is or there is not a line below it.
If the point exists we carry on. If it does not exist we mark the parameter for the line. So here there is a bunch of lines here, and it turns out -- this is Norton's straight theorem -- that there are 12 values for which you can't reduce them further. 
These are
\[
23,\ 19,\ 17,\ 11,\ 9,\ 8,\ 7,\ 5,\ 4,\ 3,\ 2,\ 1 \ . 
\]
These 12 values are the values of $k$ so that every coefficient of a replicable function is a polynomial in these 12 values of $h_k$. This is called {\em Norton basis}, and that's a very fundamental result. 
Now, Conway was talking yesterday about the work of Atkin, Fong and Smith \cite{AFS}, and by the way, Borcherds. 
Well, if Atkin, Fong and Smith had got this theorem at the time they did their computations they would have only needed at worst 24 coefficients to establish the result of the moonshine conjecture and these modular functions. 
But they didn't have it at the time and so that wasn't accessible to them. This is the general version, and this provides a basis for all replicable functions. 

I want to say a bit more about this, something special, in a minute. If the function has odd level, that means for the moonshine functions and that means that for the conjugacy class containing the element has odd order, you don't need more than these few at the bottom, in fact, $1,\ 2,\ 3,\ 5$ are sufficient to do the thing. 
There is something \cite{CohMcK} 
rather special when you have odd level, but in general you need all the above 12 elements. I'll show a quite neat proof of the theorem in a minute. This is a restatement of the condition of replicability by Norton saying that
\begin{equation}
h_{m,n}=h_{\rm{lcm}(m,\ n),\gcd(m,n)} \ .
\end{equation}
Think of the Smith normal form of a $2 \times 2$ matrix perhaps. 
You have a matrix with $m$ and $n$ on the diagonal and it's equivalent to a matrix with $\gcd(m,n)$ and $\rm{lcm}(m,\ n)$ on the diagonal. 
If you check the coefficients $h_i^{(k)}$, where $f^{(k)}(z)=\sum h_i^{(k)}q^i$, you can see that
\begin{equation}
h_i^{(k)}=k\sum_{d|k}\mu(d)h_{k/d,dki},\quad i>0,\ h_{-1}^{(k)}=1,\ h_0^{(k)}=0\ .
\end{equation}
From this, inverting the above, you deduce for all $r,s\in\BN$
\begin{equation}
h_{r,rs}=\sum_{d|r}\frac{1}{d}h_{r^2s/d^2}^{(d)}\ .
\end{equation}
This can be rewritten as
\begin{equation}
h_{m,n}=\sum_{d|(m,n)}\frac{1}{d}h_{mn/d^2}^{(d)}\ .
\end{equation}
This is the final result. 

Now those of you who have read the useful book \emph{A course in Arithmetic} by Serre \cite{Serre}, the second half of the book is devoted to things of interest to us (modular forms), and mentions the Leech lattice and various things to do with theta functions, you will find a formula very like this without the superscript $(d)$, and if you follow Serre's proofs they go through pretty well word for word in this more general situation. 

I'd like to just show you the proof of the Norton basis theorem \cite{Norton} done by Cummins \cite{Cummins} because it's very neat, and maybe other versions around are not as neat as this \footnote{
Cf.~\cite{McSev} for more discussions on the algorithms.
Of course, one sees the beginning of the supersingular primes here.
}.
\begin{theorem}
 [Norton Basis Theorem]
The $q$-coefficients of a replicable function are polynomials in
$h_k,\ k\in B=\{1,\ 2,\ 3,\ 4,\ 5,\ 7,\ 8,\ 9,\ 11,\ 17,\ 19,\ 23\}$. 
\end{theorem}
So what we want to do is to prove that for everything that is not in the basis we can actually reduce the sum $m+n$. So you need to prove that there exists $m,\ n,\ m',\ n'\in\BN$ such that
\begin{itemize}
\item[1.] $m+n=N$, where $N$ is given;
\item[2.] $\lcm(m, n) =\lcm(m',\ n')$;
\item[3.] $\gcd(m,n)=\gcd(m',n')$;
\item[4.] $m'+n'<m+n$.
\end{itemize}
\begin{proof}
The first remark is that if we find some result which is true for any number $N$ then the result can be true for $kN$. That's a useful thing to look at, and that means that we don't need to have common factors in the subscripts. Here are the cases to go through one by one:
\begin{itemize}
\item[i.]  $N=2^k$ not 2, 4 or 8; we can always look at what we believe to be the basis and see what we need not worry about. For $N=16$, here is a pair $m=1,\ n=15,\ m'=3,\ n'=5$. All else is a power of $2=16t$.
\item [ii.] For $N$ odd and $2^a+1,\ a\geq4,\ N\geq17$, we have $m=2^a-2,\ n=3,\ m'=2^{a-1}-1,\ n'=6$.
\item[iii.] $N$ odd $N\neq1+2^k$, you subtract 1 from it to get $m=N-1,\ n=1,\ m'2^{-r}(N-1),\ n'=2^r$ with $2^r | (N-1)$ is a whole divisor of $N-1$. You can follow that through $N-1>2^r\ \Rightarrow\ N>2^r+1\ \Rightarrow\ N(2^r-1)>2^{2r}-1\ \Rightarrow\ N>2^{-r}(N-1)+2^r$.
\item[iv.] Then $N$ even not a power of 2. In this situation $N$ is a going to be a product of 2, 4 or 8 with 3, 5 or 9. 
And these cases we look at individually. 
For 40 take $m=1,\ n=39,\ m'=3,\ n'=13$. For 36, $m=1,\ n=35,\ m'=5,\ n'=7$ also works for 72.
\end{itemize}
This reduces all the ones that can be reduced, and what you have left over is in the Norton basis, and this proves what the Norton basis actually is.
\end{proof}

\subsection{Elastica}
Now something of interest that maybe someone throw some light on here. If you write down the Norton basis in this way

\begin{tabular}{|p{0.2cm}|p{0.2cm}|p{0.2cm}|p{0.2cm}|p{0.2cm}|p{0.2cm}|p{0.2cm}|p{0.2cm}|p{0.2cm}|p{0.2cm}|p{0.2cm}|p{0.2cm}|p{0.2cm}|p{0.2cm}|p{0.2cm}|p{0.2cm}|p{0.2cm}|p{0.2cm}|p{0.2cm}|p{0.2cm}|p{0.2cm}|p{0.2cm}|p{0.2cm}|p{0.2cm}|}
  \hline
  & 1 & 2 & 3 & 4 & 5 &  & 7 & 8 & 9 & & 11 &  &  & & &  & 17 &  & 19 &  &  &  &23  \\
  \hline
%
  \hline
  23& &  &  & 19 &  & 17 &  &  & & &  &11 &  &9 &8 & 7 &  & 5 & 4 & 3 &2  &1  &  \\
  \hline
\end{tabular}

we have a symmetry \footnote{
In that a number is accompanied by a gap, and vice versa, except for positions 4 and 10.
} in the Norton basis with the exception of the
boxes that are supposed to contain 4 and 10. This was noticed by
Matsutani who is in Yokohama, and he wondered whether this has
to do with Weierstrass gaps \cite{Mats3}. Some of you might know about Weierstrass gaps and the implications of it if there is some connection but that I don't know \footnote{
We recall the statement of the Weierstra\ss\ gap theorem:
For a compact genus $g>0$ Riemann surface $X \ni x$, there exists exactly $g$ numbers $1 = n_1 < n_2 < \ldots < n_g < 2g$ such that there does not exist a holomorphic function on $X \backslash x$ with a pole of order $n_i$ at $x$.
}.
A remark to make is that when you are working with odd level functions, these 1, 2, 3 and 5 are the relevant entries that you need for the basis. 
Now what about 4 and 10? 
Well, you can write down expressions for 4 and 10. 
Basically speaking you can express $a_{10}$ in terms of $a_4$, and I think there might well be an argument for putting 10 in the above tables and leaving things as they are otherwise, as you can express $a_4$ in terms of $a_{10}$. 
There is a little complication which I don't want to talk about. 

This is an interesting remark and this leads to some other work of Matsutani \cite{Matsutani} which suggests that there might be some connection with the variational problem of Euler and these replicable functions. The variational problem of Euler is what you get when you take an \emph{elastica}, and elastica is Euler's word.
It's about 1750 or so. 

What is an elastica? 
Take a metal ruler and push it in from the ends it will bend, and the question is what is the curve that you get when you bend it and that is a variational problem on the integral on the square of the curvature over the arc length. 
That was solved completely by Euler \footnote{
q.v.~
Leonhard Euler, ``Methodus inveniendi lineas curvas maximi minimive proprietate gaudentes, sive solutio problematis isoperimetrici lattissimo sensu accepti, chapter Additamentum 1'', eulerarchive.org, E065, 1744.
}. 
But, if you generalized it a bit into what Matsutani calls a quantized version \cite{Matsutani2} you get some interesting objects, and genus 0 functions come up, rather than replicable functions. 
But whether there is a connection or not, I don't know
\footnote{
The reader is referred to \cite{CG} for a discussion on the significance of the genus zero property.
}
. 
So that's a curiosity which might be worth pursuing. 

\subsection{Faber Polynomials and Symmetric Functions}
Let me come back to the Faber polynomials. They haven't been much studied, really. I'd like to mention some things here. I think the first of the following identities is perhaps the most important one to remember because it's easy to remember
\begin{equation}\label{eqn:gen-hom-sym-fcns}
1+\sum_{n\geq1}h_nt^n=\prod_{i}(1-x_it)^{-1}=\exp\left(\sum
\frac{t^n}{n}\cdot\sum_i x^n_i \right)\ .
\end{equation}
The $h_n$'s are called the complete homogeneous symmetric functions \footnote{
It is also called the plethystic exponential and has been the key to a programme of counting gauge-invariant operators in quantum field theories \cite{Feng:2007ur}.
}.  
Symmetric functions can be expressed in terms of $x_1,\ x_2\cdots$, and the $h_n$ is of degree $n$ and is a sum over the $x_i$'s whether or not they are equal. 
So you are looking at sums of $x_i$'s where $x_{i_1}\geq x_{i_2}\geq\cdots$. The functions on the right hand side of \eqref{eqn:gen-hom-sym-fcns} are generating functions for these homogeneous symmetric functions \footnote{
This has been interpreted as fugacity-inserted plethystic exponential of a Hilbert series in the context of D-brane gauge theories \cite{Benvenuti:2006qr} and as Witt vectors in \cite{FM}.
}. 

McDonald \cite{McD} on symmetric functions shows that there is an involution that  transforms the homogeneous symmetric functions to the elementary ones, and thus involution changes the generating functions. 
So there is a relation between
\begin{equation}\label{eqn:gen-elem-sym-fcns}
1+\sum_{n\geq1}b_nt^n=\prod_{i}(1+x_it)=\exp\left(\sum
\frac{-t^n}{n}\cdot\sum_i x^n_i \right)
\end{equation}
and \eqref{eqn:gen-hom-sym-fcns}.

Now, what are the Faber polynomial doing?
Well, there is a completely different notion from what I've been talking about and that is that the Faber polynomials are related to a change of basis for a symmetric function. 
There are six standard bases for symmetric functions, five of which are well known, one of which is the Doubilet basis, and is called the forgotten symmetric functions \cite{Dou}.
Anyhow, one could take one of the above functions and multiply it by $1/t$ so that everything is shifted by 1. You take the following matrix (you have to be careful about these Faber polynomials)
\begin{equation}
A_n={ \left(
    \begin{array}{ccccc}
      b_1 & 1 & 0 & \cdots &0\\
      2b_2 & b_1 & 1 & \ddots &\vdots\\
       3b_3 & b_2 & \ddots & \ddots& 0\\
      \vdots & \vdots & &\ddots & 1 \\
      nb_n & b_{n-1} & && b_1 \\
    \end{array}
  \right)}\ .
\end{equation}
We have $F_n(b_1,\cdots,b_n)=\det (A_n)$  \footnote{
That's the matrix, and I got into trouble when talking to Serre about this because we were using different notations and he objected very strongly to this.
Note that the notation in footnote \ref{ft:Faber} is the unshifted version. 
}. 
Write $F_n(z)=F_n(z,b_2,\cdots,b_n)$ with $F_0(f)=1,\ F_1(f)=f,\ F_2(f)=f^2-2b_2,\ F_3(f)=f^3-3b_2f-3b_3$. 
One can think of the Faber polynomial as $F_n(z)$, and remember that I shifted everything by 1 in \eqref{eqn:gen-elem-sym-fcns} (by dividing by $q$); the function $f$ we are dealing with here is
\[
f(q)=\frac{1}{q}\left(1+\sum_{n\geq1}b_nq^n\right)\ .
\]
For replicable functions $b_1=0$ and $b_k=h_{k-1}$.
Moreover, the $F_n(z)$ are isobaric, meaning homogeneous in the subscripts, in the sense that if you replace $z$ by $b_1$ these are isobaric polynomials in the $b_i$'s. It's very easy to get signs wrong, for me anyhow, and if you get things correct for the third one I think you are ok. These are Faber polynomials, and they come about from solving the Newton relations which I just described in terms of the roots of the polynomial $F_n(z)$ and its coefficients \footnote{
That is, the recursion relations for the Faber polynomials described in footnote \ref{ft:Fn}.
}. So that's what these Faber polynomials really are.

As I said yesterday, historically it's quite interesting that they were described by Faber in 1903 and mathematicians know them, and there are certainly due to someone earlier than Faber, but it might be one of this fairly folklorist things that goes back and that might be even predates Newton. 
There isalso  Girard \footnote{
Albert Girard (1595-1632), worked on fundamental theorem of algebra, symmetric polynomials, Fibonacci numbers, inter alia.
}, but I don't know quite what role he had to play. Anyhow, Fa\`a di Bruno is the person who predates Faber, and I don't know what date we are talking about, probably $1857$; he's well known for the $n$-th derivative of the composition of two functions \footnote{
Francesco Fa\`{a} di Bruno (1825-88), cf.~
``Sullo sviluppo delle Funzioni'', Annali di Scienze Matematiche e Fisiche, 6: 479-480, 1855 and ``Note sur une nouvelle formule de calcul differentiel'', The Quarterly Journal of Pure and Applied Mathematics, 1: 359-360, 1857.
}.

\subsection{Norton's Conjecture on Replicable Functions}
Now, I'd like to make a statement on this main outstanding conjecture about replicable functions and nobody, as far as I know, has tried to solve it, but it's not quite as simple as one might wish. This is Norton's main conjecture and it is that:
\begin{conjecture} [Norton's Conjecture]
A function $f$ of the form \eqref{eqn:f} is replicable if either
\begin{itemize}
\item[0.] $f(q)=1/q+cq$ with $c\in\{0,\ -1,\ 1\}$ as we are working with integer coefficients. Respectively for these values of $c$ we get the $\exp,\ \sin$ and $\cos$ functions \footnote{Recall that the nome $q = \exp(2 \pi i \tau)$.}, which we call ``Modular Fictions'' and to be ignored henceforth.

Or, surprisingly using the modular polynomial you can get some results about these things, which are consequences of the modular polynomial for them. 
Things like the $\cos$ of twice the angle is a polynomial of the $\cos$ whereas the $\sin$ of twice the angle is not a polynomial in the $\sin$. 
This was proved by C.~Cummins again \cite{Cummins}, that this is all there is \footnote{Indeed, these are the Chebyshev polynomials; the reader is referred to \cite{McK} for discussions on how these famous polynomials are the simplest replicable functions.}.

\item[1.] There exists $N:\ \G_0(N)\subset G_f\subset \mbox{Nor}(\G_0(N))$, where Nor is the normalizer inside  $\PR$,
such that the compact Riemann surface $\widehat{G_f \slash \fH}$ obtained by adding a finite set of inequivalent cusps has genus 0, $f$ is a principal modulus of this Riemann surface, and $G_f$ is commensurable \footnote{
Recall that $G_f$ was defined in \eqref{Gf} as the modular subgroup for which $f$ is invariant (weight $k=0$). 
} with $\PS$. 

The maximal groups on the right  side are called \emph{Helling groups}, and there is a paper by Conway \cite{conway} called ``Understanding groups like $\G_0(N)$'' which was referred to by J.~Duncan. It's really a pretty piece of work in that he also proves Helling's theorem. 
\end{itemize}
\end{conjecture}

What is needed is a proof of this result. 
I think the proper way to do this is to use a two-sided decomposition with respect to GL$_2$ of the ad\`eles, and in that way one should be able to pick up the primes dividing the Monster's order, which are the primes you find inside the above levels $N$, namely the 15 supersingular primes.

What has been done to now, I might say something about computations if there is time, it's been done by some people; I had some visitors that did some work, C.~ Cummins did some, but the bulk of this stuff was done by S.~Norton and I could describe it. Basically speaking, it was a very local picture that was used to do the computations to find these replicable functions. 
There are one or two functions that you know classically from Weber's \cite{weber} and Schl\"afli's \cite{schlafli} work, and from there you can build up other functions which are closely related to replicable functions. You just keep looking and making sure that the genus is zero. But it is conceivable that something has been left out, it's unlikely but it's possible. 

The number 616 is the number of replicable functions there are (with integer coefficients). That in itself is a quite interesting number, and there is amusing reference to J.~Conway's remark. 
Remember that he was talking about\footnote{
The group is defined as follows.
Consider a graph $G_{p,q,r}$ with a single tri-vertex, say $a$ and 3 strands, each consisting of respectively $p$, $q$ and $r$ nodes joined up; call these $b_{1,\ldots,p}$, $c_{1,\ldots,q}$ and $d_{1,\ldots,r}$.
Thus there are $p+q+r+1$ nodes in total.
The group $Y_{p,q,r}$ is a Coxeter-type group with one generator associated to each of the nodes and presentation
\[
Y_{p,q,r} = \gen{a, b_{1,\ldots,p}, c_{1,\ldots,q}, d_{1,\ldots,r} |
(gg')^{O(g,g')} = (ab_1b_2ac_1c_2ad_1d_2)^{10} = 1
} \ ,
\]
where $g$ and $g'$ are any of the $a,b,c,d$ generators and that $O(g,g')=3$ (respectively 2) if $g$ and $g'$ are adjacent (respectively, not adjacent).
It was shown \cite{ivanov} that $Y_{5,5,5} \simeq Y_{4,4,4} \simeq \IM \wr C_2$, the wreath product of the monster with the cyclic group of order 2.
} a group $Y_{555}$ or perhaps $Y_{666}$. Well, $666$ as you know occurs in the Bible, it's sort of a bad number, and I was giving this talk in Norway to the mathematical department and I mentioned this $616$ as the number of replicable functions and a person in the back went out and he came back to me with a transparency on which there is something called the ``Oxyrhynchus Papyrus'' from few hundred A.~D.~\footnote{
Manuscripts discovered in the C19th near Oxyrhynchus, Egypt, dating from 1st to 6th century AD., mostly housed in the Ashmolean in Oxford.
}
In the papyrus, which is written in Greek-- it's very readable, it says that the real number shouldn't be 666; it should be 616. So that was very quite amusing. I don't know why or how he came across to this.

Here is another piece of numerology, entirely frivolous but nevertheless grounded on some deep mathematics, which may amuse you.
The Leech lattice $\Lambda_{24}$, the famous even unimodular lattice in 24 dimensions, can be constructed from the even Lorentzian unimodular lattice $II_{25,1}$ in 26 dimensions using a Weyl vector $w = (0,1,2,3,\ldots,23,24;70)$.
Then $\Lambda_{24}$ is realized as $w^\perp / w$.
That $w$ is indeed an integral vector in $II_{25,1}$ follows from the remarkable Diophantine condition \footnote{
Discussions on this equation and the emergence of 42 from $j(q)$ are presented in \cite{He:2014uma}, in a volume in honour of J.~H.~Conway. The number of pages of the present notes is, of course, 42.
}
\begin{equation}
1^2 + 2^2 + 3^2 + \ldots + 23^2 + 24^2 = 70^2 \ .
\end{equation}
Now, consider the first 24 $q$-series coefficients of the (normalized) modular $J$-function, viz, 
$196884, 21493760, 864299970, 20245856256, \ldots$ (that is we do not include the constant term 744 nor $1/q$ and start from $\cO(q)$).
Sum the squares of these 24 numbers and compute it modulo 70, you will get 42, which we all know to be the answer to the ultimate question of life, the universe and everything.
So, here again, you see how Moonshine encodes all things.

\subsection{Mahler's Recurrence Relations}\label{s:mahler}
Now, there is another interesting paper \cite{ACMS} where we used recurrences to build up these series. 
What is not very well known is the recurrence that we use which people call Borcherds' recurrence \footnote{
The sort of recurrences which arise from Borcherds' proof is the remarkable one such as
\[
j(p) - j(q) = \left( \frac{1}{p} - \frac{1}{q} \right)
\prod\limits_{m,n=1}^\infty (1 - p^n q^m)^{c_{nm}} \ ,
\]
where $c_k$ are the q-expansion coefficients of $j(q)$.
We will see more recurrences in \S\ref{s:mahler}.
}, which is fair enough, is really due to Mahler and Mahler wrote his paper in 1974, published in 1976 in the journal of the Australian Math.~Soc. \cite{Mahler}.

He was a very remarkable number theorist. He was a cripple and his father knew Carl Ludwig Siegel. 
So he got to get a rather good education without going through the usual processes of university and then he became a refugee in Manchester for 30 years and then he went from an assistant professor to a personal Chair, in Canberra where he lived and died not a long ago. 

He wrote me a letter in 1982 or 1983 in which he said that he thought his paper has something to do with moonshine, and if I would come to Canberra and discuss it with him before he dies. 
I was in Montreal and he was in Australia. I went out there and unfortunately I got ill and I came back. So we never actually met, but Mahler recurrence is the recurrence that does the calculations.

We observe that (q.v.~also \cite{He:2015yoa})
\begin{align}
\nn
360+256 = &616 \\
120+2\times248= &616 \ .
\end{align}
Maybe you can explain this. 
I'll give you the Mahler recurrence.  I'm just going back a bit. 
I already remarked upon
\[ 
f^{(n)}(nz)=\frac{1}{n}F_n(f(z))-\frac{1}{n}
     \sum_{\substack{ad=n\\0\leq b<d\\a<n}}
     f^{(a)}\left(
     \frac{az+b}{d}
     \right)\ ,
\]
and that $f^{(n)}(nz)$ is invariant under $z\ \mapsto\ z+1/n$.
Cummins and Norton have proved the replicability of rational Hauptmoduln \cite{Cum-Nor}, and remember that the replicability property is essentially combinatorial. You are just saying that things work nicely with Hecke operators whereas if you want to prove results about replicable being Hauptmodul you have to go from the combinatorial side to the analytic and it's a much more difficult thing to do. 

There is some very remarkable work done by a student of Arne Meurman in Sweden. His name is Dmitry Kozlov \cite{Kozlov}. Masao Koike is the first to realize what the generalized Hecke operator is. And if I have done it correctly. When $n=p$ is prime, the classical Hecke operator can be expressed in terms of the the Atkin $U_-$-operator and if you know it for $p$ you know it for all integers because of the multiplicative property, and the Adams' $V_-$-operator
\begin{align}
\nn
T_p&=\frac{1}{p}V_p+U_p \\
\nn
V_p&:\ f(q)\ \mapsto\ f(q^p)\\
U_p&:\ a_nq^n\ \mapsto\ a_{pn}q^n\ .
\end{align}
Koike recognizes this for the twisted Hecke operator in the very early days
\begin{equation}
\widehat{T}_p=\frac{1}{p}\Psi^p\circ{V}_p+U_p
\end{equation}
with $\widehat{V}_p=\Psi^p\circ V_p:\ f(q)\ \mapsto\ f^{(p)}(q^p)$. 
Now if you go back to the moonshine to remember what this does for you, this means that the coefficients are traces which are the sums of the initial $p-$th power function. 
That's the context of moonshine but the rest is true for replicable functions generally.

There are some characteristic classes associated with this, and this is called the {\em Bott Cannibalistic Class}.  
So when I was in Harvard I asked Bott about it, but he told me he had forgotten so I don't think we'll get any further there. 
By the way there is a remarkable meeting in a week or two in Montreal, I think in June, a week on Bott's legacy with very good speakers including Witten and Atiyah among others in the University of Montreal \footnote{
This is the conference ``A Celebration of the Mathematical Legacy of Raoul Bott'', CRM, June 9-13, 2008, Montreal.
One can find the proceedings in {\it CRM Proceedings \& Lecture Notes}, Volume 50, AMS 2010.
}. 
I just make a passing remark is that
$f\equiv f^{(p)}\bmod\,p$ because in $f^{(p)}$ the coefficients are sums of $p$-th powers of the coefficients of $f$. This is the work of Koike, and there is a nice survey article by him in Japanese in Sugaku, and there is an English version in number 160 of the AMS translations \cite{koike}.

I think Mahler must be the first or the only mathematician, unless Conway has done the same thing, to publish a calculator program in the Royal Society of London Proceedings, but he did that in an subsequent paper where he does his computations. 
What's extremely interesting about Mahler is that for each prime $p$ he has a recurrence relation to compute the coefficients of these functions and he is really computing the $q$-coefficients of the $j$-function. 
There are several people who have done this for $j$, where the actual setup will work for replicable functions in general. This was true by Kozlov's thesis which was on the $j$-function, but in fact it applies to all replicable functions. The same with Mahler; he has this recurrence relations for all prime $p$ for the $j$-function and he tries to use them in a broader context to other functions that have arisen in \cite{fricke} and others and it doesn't work and he doesn't see why it doesn't work. 

Today we know why it doesn't work. It's because the level of the function and the level of the Hecke operator are not coprime, and that's the modification that you need to go from Borcherds' formula to Mahler's formula. 
So I'm just going to give it to you here. 
For $p=2$ Mahler has used
\begin{align}
\nn
f\left(\frac{z}{2}\right)+f\left(\frac{z+1}{2}\right)+f(2z)&=f(z)^2-2a_1
\\
f\left(\frac{z}{2}\right)f(2z)+f\left(\frac{z+1}{2}\right)f(2z)+
f\left(\frac{z}{2}\right)f\left(\frac{z+1}{2}\right)&=2a_2f-f+2(a_4-a_1)\ ,
\end{align}
which might be look at as the elementary symmetric functions of degree 1 and 2 respectively in the function values $f\left(\frac{z}{2}\right),\ f\left(\frac{z+1}{2}\right)$ and $f(2z)$. This is for the $j$-function originally, and then you have to make modifications. The modifications are simple to make; whenever you see a 2 in the $f$-value, say in $f(2z)$ you change it to $f^{(2)}(2z)$. These are the recurrence relations from Mahler
\comment{
\begin{align}
\hspace{-2cm}
\nn
a_{4k}&=\sum_{j=1}^{k-1}a_ja_{2k-j}+\frac{1}{2}(a_k^2-{a_k}_\bullet) \ , \\
\nn
a_{4k+1}&=a_{2k+3}\sum_{j=1}^{k}a_ja_{2k+2-j}+\frac{1}{2}
(a_{k+1}^2-{a_{k+1}}_\bullet)+\frac{1}{2}(a_{2k}^2-{a_{2k}}_\bullet)-a_2a_{2k}+\\
\nn
& 
\qquad \qquad
+ \sum_{j=1}^{k}{a_j}_\bullet a_{4k-4j}+\sum_{j=1}^{2k-1}(-1)^ja_ja_{4k-j} \ , \\
\nn
a_{4k+2}&=a_{2k+2}+\sum_{j=1}^{k-1}a_ja_{2k+1-j} \ , \\
a_{4k+3}&=a_{2k+4}\sum_{j=1}^{k+1}a_ja_{2k+3-j}+\frac{1}{2}
(a_{2k+1}^2-{a_{2k+1}})-a_2a_{2k+1}+
\sum_{j=1}^{k}{a_j}_\bullet a_{4k+2-4j}+\sum_{j=1}^{2k}(-1)^ja_ja_{4k+2-j}
\end{align}

The dot 
means to put the superscript $(2)$ there. So the recurrence relation for replicable functions are
}
\begin{align}
\nn
a_{4k}&=\sum_{j=1}^{k-1}a_ja_{2k-j}+\frac{1}{2}(a_k^2-a_k^{(2)}) \ ,\\
\nn
a_{4k+}&=a_{2k+3}\sum_{j=1}^{k}a_ja_{2k+2-j}+\frac{1}{2}
(a_{k+1}^2-a_{k+1}^{(2)})+\frac{1}{2}(a_{2k}^2-{a_{2k}^{(2)}})-a_2a_{2k}+\\
\nn
& \qquad \qquad 
+ \sum_{j=1}^{k}{a_j^{(2)}}a_{4k-4j}+\sum_{j=1}^{2k-1}(-1)^ja_ja_{4k-j} \ ,\\
\nn
a_{4k+2}&=a_{2k+2}+\sum_{j=1}^{k-1}a_ja_{2k+1-j}\ , \\ 
a_{4k+3}&=a_{2k+4}\sum_{j=1}^{k+1}a_ja_{2k+3-j}+\frac{1}{2}
(a_{2k+1}^2-{a_{2k+1}})-a_2a_{2k+1}+ \\
\nn
& \qquad \qquad 
+ \sum_{j=1}^{k}{a_j^{(2)}}a_{4k+2-4j}+\sum_{j=1}^{2k}(-1)^ja_ja_{4k+2-j} \ .
\end{align}

which do come from
\begin{align}
\nn
f\left(\frac{z}{2}\right)+f\left(\frac{z+1}{2}\right)+f^{(2)}(2z)&=f(z)^2-2a_1
\\
f\left(\frac{z}{2}\right)f^{(2)}(2z)+f\left(\frac{z+1}{2}\right)f^{(2)}(2z)+
f\left(\frac{z}{2}\right)f\left(\frac{z+1}{2}\right)&=2a_2f-f+2(a_4-a_1)\ ,
\end{align}
and these are Borcherds' relations. That's the difference between Mahler's and Borcherds'. 

Notice that the original ones are universal, and if you have any function in this case of odd level
these are the recurrence relations for its coefficients and if you have a function of even level the you have to know what $f^{(2)}$ is in order to make use them; you have to use coefficients from $f^{(2)}$ in order to build up these relations. They come in a  group of four. They are extremely good for computing several hundreds coefficients quite easily, but if you speak to someone like Atkin few hundreds is nothing; he computes the first 10,000.

Now, to go back to some remarks a bit earlier on. 
We find $f^{(p)}=\sum_{n}h_n^{(p)}q^n$ and
\begin{equation}
h_n^{(p)}=ph_{pn,p}-pa_{p^2n}\ .
\end{equation}
On $\mathbb{M}:\ f_{\gen{g}}\ \longrightarrow\ f_{\gen{g}}^{(k)}=f_{\gen{g^k}}$, this is the replication formula I mentioned to you \footnote{
Thus we come full circle back to the Monster $\IM$.
Recall that for conjugacy classes, and in particular of cyclic groups $\gen{g}$, one can write down McKay-Thompson series, these obey replication identities.
The question then is whether group structure in $\IM$, such as the power map $g \mapsto g^k$ are reflected in the replicable functions $f\mapsto f^{(k)}$.
}. Since $|\mathbb{M}|<\infty$ there exists $k=k_0$ such that $g^{k_0}=id$ and so $f_{\gen{g^{k_0}}}=f_{\gen{id}}=j$. So in a sense you can think of the functions we start with as replication roots of the $j$-function.


\section*{Acknowledgements}
JM is grateful to the organizers for the Kashiwa conference for their warm hospitality and that of the NSERC of Canada.
YHH would like to thank the Science and Technology Facilities Council, UK, for grant ST/J00037X/1, and the Chinese Ministry of Education, for a Chang-Jiang Chair Professorship at NanKai University.



\begin{thebibliography}{aaaa}

\bibitem[ACMS]{ACMS}
D.~Alexander, C.~Cummins, J.~McKay, C.~Simons, 
``Completely replicable functions,'' in Liebeck, Saxl, 
``Groups, Combinatorics and Geometry'', 
LMS Lecture Note Series (CUP) 165: 87–98.


\bibitem[AFS]{AFS}
S.~D.~Smith, ``On the head characters of the Monster simple group'',
Finite Groups – Coming of Age (Montr\'eal, 1982), Contemp. Math. 45
(American Mathematical Society, Providence 1996).

\bibitem[Apo]{Apo}
Tom M.~Apostol, ``Modular functions and Dirichlet series in number theory'', Second edition. Graduate Texts in Mathematics, 41. Springer-Verlag, New York, 1990.

\bibitem[AKM]{Aspinwall:2000kf} 
  P.~S.~Aspinwall, S.~H.~Katz and D.~R.~Morrison,
  ``Lie groups, Calabi-Yau threefolds, and F theory,''
  Adv.\ Theor.\ Math.\ Phys.\  {\bf 4}, 95 (2000)
  \href{http://arxiv.org/abs/hep-th/0002012}{[hep-th/0002012]}.


\bibitem[Atk-Swi]{Atk-Swi} 
A.O.L.~Atkin \& H.P.F.~Swinnerton-Dyer, `` Modular forms on noncongruence subgroups'', Combinatorics (Proc. Sympos. Pure Math., Vol. XIX, Univ. California, Los Angeles, Calif., 1968), pp. 1-25. AMS., Providence, R.I., 1971.

\bibitem[ATLAS]{atlas}
J.H. Conway, R.T. Curtis, S.P. Norton, R.A. Parker, R.A. Wilson,
``Atlas of Finite groups'', 
Oxford University Press, Eynsham, UK, 1985.
\href{http://brauer.maths.qmul.ac.uk/Atlas/v3/}{http://brauer.maths.qmul.ac.uk/Atlas/v3/}

\bibitem[BFHH]{Benvenuti:2006qr} 
  S.~Benvenuti, B.~Feng, A.~Hanany and Y.~-H.~He,
  ``Counting BPS Operators in Gauge Theories: Quivers, Syzygies and Plethystics,''
  JHEP {\bf 0711}, 050 (2007)
\href{http://arxiv.org/abs/hep-th/0608050}{[hep-th/0608050]}


\bibitem[Blog$^>$]{blogs}
Some blogs on Moonshine: \\
\url{http://ncatlab.org/nlab/show/Moonshine} \\
\url{http://www.neverendingbooks.org/index.php/monsters-and-moonshine-a-booklet.html}


\bibitem[Bor]{borcherds} 
R.~E.~Borcherds,
  ``Vertex algebras, Kac-Moody algebras, and the monster,''
  Proc.\ Nat.\ Acad.\ Sci.\  {\bf 83}, 3068 (1986);
  ``Monstrous moonshine and monstrous Lie superalgebras'', 
Invent.~Math.~109 (1992) 405–444.

\bibitem[Cam]{cam}
P.~J.~Cameron, {\it Permutation Groups}, LMS Student Texts, 45, CUP (1999). 

\bibitem[Chev]{Chev} C.~Chevalley, ``Sur certains groupes simples''. 
(French) Tohoku Math. J. (2) 7 (1955), 14-66.

\bibitem[CohMck]{CohMcK}
Harvey Cohn, John McKay,
``Spontaneous Generation of Modular Invariants''
Math. Comp. 65 (1996), 1295-1309.


\bibitem[Con]{conway}
J.H.~Conway, ``Understanding groups like $\Gamma_0(N)$''. 
Groups, difference sets, and the Monster (Columbus,. OH, 1993), 327-343.

\bibitem[Con-Nor]{Con-Nor} 
J.H.~Conway \& S.P.~Norton, ``Monstrous moonshine'', Bull. London Math. Soc. 11 (1979), no. 3, 308-339.

\bibitem[Cum]{Cummins}
C.~J.~Cummins, ``Some comments on replicable functions'', Modern trends in
Lie algebra representation theory (Queen’s Univ., Kingston, ON, 1994) 48–55,
Queen’s Papers in Pure and Appl. Math. 94 (1994).

\bibitem[CumGan]{CG}
C.~J.~Cummins and T.~Gannon, ``Modular equations and the genus zero property of moonshine functions,'' Inventiones mathematicae, 1997, Volume 129, Issue 3.

\bibitem[Cum-Nor]{Cum-Nor}
C.J.Cummins and S.P.Norton, ``Rational Hauptmoduls are replicable'', 
Canad. J. Math. 47 (1995), no. 6, 1201–1218.

\bibitem[CDH$^>$]{Cheng:2012tq} 
  M.~C.~N.~Cheng, J.~F.~R.~Duncan and J.~A.~Harvey,
  ``Umbral Moonshine,''
  \href{http://arxiv.org/abs/1204.2779}{arXiv:1204.2779 [math.RT]}.

\bibitem[CDHKW$^>$]{Cheng:2013kpa} 
  M.~C.~N.~Cheng, X.~Dong, J.~Duncan, J.~Harvey, S.~Kachru and T.~Wrase,
  ``Mathieu Moonshine and N=2 String Compactifications,''
  \href{http://arxiv.org/abs/1306.4981}{arXiv:1306.4981 [hep-th]}.

\bibitem[Che$^>$]{Cheng:2010pq} 
  M.~C.~N.~Cheng,
  ``K3 Surfaces, N=4 Dyons, and the Mathieu Group M24,''
  Commun.\ Num.\ Theor.\ Phys.\  {\bf 4}, 623 (2010)
  \href{http://arxiv.org/abs/1005.5415}{[arXiv:1005.5415 [hep-th]]}.


\bibitem[de Bra]{de bran} 
L.~de Branges, ``Underlying concepts in the proof of the Bieberbach conjecture'', A plenary address presented at the International Congress of Mathematicians held in Berkeley, California, August 1986. Introduced by Max M. Schiffer. ICM Series. American Mathematical Society, Providence, RI, 1988.


\bibitem[Deg-Wen$^>$]{Wend} 
A.~Degeratu \& K.~Wendland, ``Friendly giant meets pointlike instantons? On a new conjecture by John McKay'', 
Moonshine: the first quarter century and beyond, 55-127, London Math. Soc. Lecture Note Ser., 372, Cambridge Univ. Press, Cambridge, 2010.

\bibitem[DGO$^>$]{DGO}
John F.~R.~Duncan, Michael J.~Griffin, Ken Ono, ``Moonshine'', Res.~ in the Math.~Sciences (2015) 2:11,
arXiv:1411.6571 [math.RT]

\bibitem[DGO2$^>$]{DGO2}
John F.~R.~Duncan, Michael J.~Griffin, Ken Ono, ``Proof of the Umbral Moonshine Conjecture'',
Res.~ in the Math.~Sciences (2015) 2:26, 	arXiv:1503.01472 [math.RT]

\bibitem[DKM]{DKM}
D.~Dummit, H.~Kisilevsky, and J.~McKay, ``Multiplicative products of eta functions'', in {\it Finite groups - coming of age} (Montreal, Que., 1982),
vol. 45 of Contemp. Math., pp. 89–98. Amer. Math. Soc., Providence, RI,
1985. (The reviewer in Math Reviews points out and corrects a
printing error in the paper).

\bibitem[Dou]{Doubilet}
P.~Doubilet, ``On the Foundations of Combinatorial Theory. VII: Symmetric Functions through the theory of distribution and occupancy'', Studies in Applied Maths, Vol LI, 4, 1972.

\bibitem[DO$^>$]{do}
John F.~R.~Duncan, Ken Ono,
``The Jack Daniels Problem'', J.~Number Theory 161 (2016) pp.230 - 239,
arXiv:1411.5354 [math.NT]


\bibitem[duS$^>$]{Sautoy}
Marcus du Sautoy,
``Finding Moonshine: A Mathematician's Journey Through Symmetry'',
 Harper Perennial, 2009, ISBN-13: 978-0007214624.

\bibitem[EOT$^>$]{Eguchi:2010ej} 
  T.~Eguchi, H.~Ooguri and Y.~Tachikawa,
  ``Notes on the K3 Surface and the Mathieu group $M_{24}$,''
  Exper.\ Math.\  {\bf 20}, 91 (2011)
  \href{http://arxiv.org/abs/1004.0956}{[arXiv:1004.0956 [hep-th]]}.


\bibitem[Erd]{Erd} 
C.~Erdenberger, ``The Kodaira dimension of certain moduli spaces of abelian surfaces'', Math. Nachr. 274/275 (2004), 32-39. 

\bibitem[FHH]{Feng:2007ur}
B.~Feng, A.~Hanany and Y.~H.~He,
``Counting gauge invariants: The Plethystic program,''
JHEP \textbf{03} (2007), 090
[arXiv:hep-th/0701063 [hep-th]].

\bibitem[FM$^>$]{FM}
Roland Friedrich, John McKay
``Formal Groups, Witt vectors and Free Probability'', 
\href{http://arxiv.org/abs/1204.6522}{arXiv:1204.6522}.


\bibitem[Fra]{frame}
J.~S.~Frame, ``The Theory of Tables of Group Characteristics,''
Harvard Univ.~thesis, 1933.


\bibitem[Fri]{fricke}
R.~Fricke,  ``Die elliptischen Funktionen und ihre Anwendungen''. 
Zweiter Band, 1922; Reprint, Springer Verlag, 2011.

\bibitem[GHV$^>$]{Gaberdiel:2010ca} 
  M.~R.~Gaberdiel, S.~Hohenegger and R.~Volpato,
  ``Mathieu Moonshine in the elliptic genus of K3,''
  JHEP {\bf 1010}, 062 (2010)
  \href{http://arxiv.org/abs/1008.3778}{[arXiv:1008.3778 [hep-th]]}.


\bibitem[Gan]{Gannon:2004xi} 
  T.~Gannon,
  ``Monstrous moonshine: The First twenty five years,''
  \href{http://arxiv.org/abs/math/0402345}{arxiv:math/0402345 [math-qa]}.


\bibitem[GK$^>$]{Govindarajan:2009qt} 
  S.~Govindarajan and K.~Gopala Krishna,
  ``BKM Lie superalgebras from dyon spectra in Z(N) CHL orbifolds for composite N,''
  JHEP {\bf 1005}, 014 (2010)
  \href{http://arxiv.org/abs/0907.1410}{[arXiv:0907.1410 [hep-th]]}.

\bibitem[Grun]{Grun}
H. ~Grunsky, ``Koeffizientenbedingungen f\"ur schlicht abbildende meromorphe Funktionen'', Math. Z., 45 (1939) 29-61.

\bibitem[HeJe]{He:2003pq} 
  Y.~-H.~He and V.~Jejjala,
  ``Modular matrix models,''
  \href{http://arxiv.org/abs/hep-th/0307293}{hep-th/0307293}.

\bibitem[HM$^>$]{He:2012kw} 
Y.~-H.~He and J.~McKay, ``N=2 Gauge Theories: Congruence Subgroups, Coset Graphs and Modular Surfaces,'' J.\ Math.\ Phys.\  {\bf 54}, 012301 (2013)
\href{http://arxiv.org/abs/1201.3633}{[arXiv:1201.3633 [hep-th]]}.


\bibitem[HeMcK$^>$]{He:2013lha} 
  Y.~-H.~He and J.~McKay,
  ``Eta Products, BPS States and K3 Surfaces,''
\href{http://arxiv.org/abs/1308.5233}{arXiv:1308.5233 [hep-th]}.

\bibitem[HeMcK2$^>$]{He:2014uma}
Y.~H.~He and J.~McKay,
``Moonshine and the Meaning of Life,'' in {\it Contemporary Mathematics 694, Ed. M.~Bhagarva et al. 2017}
[arXiv:1408.2083 [math.NT]].

\bibitem[HeMcK3$^>$]{He:2015yoa}
Y.~H.~He and J.~McKay,
``Sporadic and Exceptional,''
[arXiv:1505.06742 [math.AG]].

\bibitem[HMR$^>$]{He:2012jn} 
  Y.~-H.~He, J.~McKay and J.~Read,
  ``Modular Subgroups, Dessins d'Enfants and Elliptic K3 Surfaces,''
 \href{http://arxiv.org/abs/1211.1931}{arXiv:1211.1931 [math.AG]}.


\bibitem[HSZ$^>$]{Haldar:2021rri}
P.~Haldar, A.~Sinha and A.~Zahed,
``Quantum field theory and the Bieberbach conjecture,''
[arXiv:2103.12108 [hep-th]].



\bibitem[Iva]{ivanov}
A.~Ivanov, ``Y-groups via transitive extension'', J.~Alg.~218 (1999) 412 - 435.

\bibitem[Kil]{kilford}
L.~J.~P.~Kilford, ``Generating spaces of modular forms with $\eta$-quotients'',
\href{http://arxiv.org/abs/math/0701478}{arXiv:math/0701478}.


\bibitem[KMP$^>$]{Mats3}
J.~Komeda, S.~Matsutani, E.~Previato,
``The sigma function for Weierstrass semigroups $\left<3,7,8\right>$ and $\left<6,13,14,15,16\right>$'',
arXiv:1303.0451 [math.AG]


\bibitem[Koi]{koike}
Masao Koike, ``Moonshine: a mysterious relationship between simple groups and automorphic functions'' in {\it Selected papers on number theory, algebraic geometry, and differential geometry}, 33–45, Amer.~Math.~Soc.~Transl. Ser. 2, 160, AMS. Providence, RI, 1994.

\bibitem[Kon]{Kon} 
Takeshi Kondo, ``The automorphism group of Leech lattice and elliptic modular functions'', J. Math. Soc. Japan 37 (1985), no. 2, 337-362.

\bibitem[Koz]{Kozlov}
Dmitry Kozlov,``On Functions Satisfying Modular Equations for Infinitely Many Primes'', Canad. J. Math. 51(1999), 1020-1034.\\
-- ``On Completely Replicable Functions and Extremal Poset Theory'',
Masters Thesis, University of Lund, 1994.

\bibitem[Mah]{Mahler}
K.~Mahler, ``On a Class of Non-Linear Functional Equations Connected with Modular Functions''. J. Aust. Math. Soc 22A (1976) 65-118.

\bibitem[Mar]{Mar} 
Y.~Martin, ``Multiplicative $\eta$-quotients''. 
Trans. Amer. Math. Soc. 348 (1996), no. 12, 4825-4856.

\bibitem[MarOno]{MarOno} 
Yves Martin and Ken Ono, ``Eta-quotients and elliptic curves'',
Proc. Amer. Math. Soc. 125 (1997)


\bibitem[Mat]{Mat} 
E.~Mathieu, ``M\'emoire sur l'\'etude des fonctions de plusieurs quantit\'es, sur la mani\'ere de les former et sur les substitutions qui les laissent invariables'', J. Math. Pures Appl. (Liouville) (2) VI, 1861, pp. 241-323.

\bibitem[Mats]{Matsutani}
Shigeki Matsutani,
``Euler's Elastica and Beyond,''
J.~Geometry and Symmetry in Physics, 17 (2010) 45-86.

\bibitem[Mats2]{Matsutani2}
Shigeki Matsutani,
``Relations in a quantized elastica'',
J.~Phys.~A. vol. 41 issue 7 February 22, 2008. p. 075201


\bibitem[McD]{McD}
Ian G.~Macdonald, ``Symmetric functions and Hall polynomials'', Second ed. Oxford Mathematical Monographs, OUP, 1995, ISBN 0-19-853489-2.



\bibitem[McK]{ade}
J.~McKay, ``Graphs, Singularities, and Finite Groups,'' Proc. Symp. Pure Math.
Vol 37, 183-186 (1980).

\bibitem[McK2]{essentials}
J.~McKay, ``The Essentials of Monstrous Moonshine'' Adv.~Studies in Pure Maths 32, 2001, pp347-353.


\bibitem[MS]{mckaysebbar}
J.~McKay and Abdellah Sebbar, ``Arithmetic Semistable Elliptic Surfaces'', 
Proceedings on Moonshine and related topics (Montr\'eal, QC, 1999), 119--130,
CRM Proc. Lecture Notes, 30, Amer. Math. Soc., Providence, RI, 2001.

\bibitem[McK-Seb]{McK} 
John McKay \& Abdellah Sebbar, ``Replicable functions: an introduction'', Frontiers in number theory, physics, and geometry. II, 373-386, Springer, Berlin, 2007.

\bibitem[McKSev]{McSev}
J.~McKay \& D.~Sevilla, ``Decomposing replicable functions''. 
LMS J. Comput. Math. 11: 146 - 171, arxiv:0803.3419[math.NT]


\bibitem[New]{New}
Morris Newman,
``Modular Forms Whose Coefficients Possess Multiplicative Properties'',
Annals of Mathematics, Vol. 70, No. 3 (Nov., 1959), pp. 478-489.

\bibitem[Mil]{Mil} 
G.A.~Miller, ``Sur plusieurs groupes simples'', 
(French) Bull. Soc. Math. France 28 (1900), 266-267.

\bibitem[Nor]{Norton}
S.~P.~Norton, ``More on moonshine'', 
Computational group theory (London Academic Press, 1984) 185–193.

\bibitem[Pomm]{Pomm} 
Ch.~Pommerenke, ``Uber die Faberschen Polynome schlichter Funktionen'', Mathematische Z., 85 (1964) 197-208.

\bibitem[PolyMath$^>$]{PolyMath}
The polymath project, \url{https://polymathprojects.org/}

\bibitem[Ogg]{Ogg}
A.~Ogg, ``Modular Functions'', In {\it The Santa Cruz Conference on Finite Groups}. Ed.~B.~Copperstein, G.~Mason, June 25–July 20, 1979. Providence, RI: Amer. M
ath. Soc. pp. 521–532.

\bibitem[Rez]{reznick}
B.~Reznick, ``Resources for Research (an always preliminary list)'',
\url{http://www.math.uiuc.edu/~reznick/rfr.html}

\bibitem[Ron]{ronan}
Mark Ronan, ``Symmetry and the Monster: One of the greatest quests of mathematics'', OUP, 2007, ISBN-13: 978-0192807236.

\bibitem[San$^>$]{sankaran}
G.~Sankaran, ``A supersingular coincidence'', 	arXiv:2009.11379 [math.NT]


\bibitem[Sch]{schlafli}
Lugwig Schl\"{a}fli, J.H.Graf ed., 
``Theorie der vielfachen Kontinuit\"{a}t'', 
Republished by Cornell University Library historical math monographs 2010 
(in German), Zürich, Basel: Georg \& Co., ISBN 978-1-4297-0481-6
(1901) [1852]

\bibitem[Seb1]{classSebbar}
A.~Sebbar, 
``Classification of torsion-free genus zero congruence groups'', 
Proc. Amer. Math. Soc. 129 (2001), 2517--2527.

\bibitem[Seb2]{sebbar}
Abdellah Sebbar, 
``Modular subgroups, forms, curves and surfaces'', 
Canad. Math. Bull. 45 (2002), no. 2, 294--308.

\bibitem[Ser]{Serre}
J-P~Serre, 
``Cours d'arithm\'etique,'' Springer 1973.



\bibitem[THM$^>$]{Tatitscheff:2018aht}
V.~Tatitscheff, Y.~H.~He and J.~McKay,
``Cusps, Congruence Groups and Monstrous Dessins,'' Indagationes Mathematicae,
31, 6, (2020), pp1015 - 1065, [arXiv:1812.11752 [math.NT]].

\bibitem[Web]{weber}
Heinrich Martin Weber, 
``Lehrbuch der Algebra'' (in German) (3rd ed.), 
New York: AMS Chelsea Publishing, ISBN 978-0-8218-2971-4; (1981) [1898].
\end{thebibliography}
\end{document}